\documentclass[11pt]{article}
\usepackage{amssymb,amsmath, amsthm}
\usepackage{graphicx}
\usepackage{graphics}
\usepackage{psfrag}

\usepackage[all]{xy}
 
\newcommand{\Rb}{\mathbb{R}}

\begin{document}
\title{Numerical Simulations in  Two-Dimensional Neural Fields}
\author{Pedro M. Lima and Evelyn Buckwar\\
Institute of Stochastics\\
Johannes Kepler University\\
Altenbergerstr. 69 , 4040 Linz, Austria}
\maketitle

\begin{abstract} 
In the present paper we are concerned with a numerical algorithm for the approximation of
the two-dimensional neural field equation with delay. This algorithm was  described in 
\cite{LP1} and \cite{LP2}. We consider three numerical examples that have been analysed before by other authors and are directly connected with real world applications.
The main purposes are 1) to test the performance of the mentioned algorithm, by comparing the numerical results with those obtained by other authors;  2) to analyse  with more detail the properties of the solutions and take conclusions about their physical meaning.
\end{abstract}
\section{Introduction}
We are concerned about the numerical solution of the following integrodifferential
equation:
\begin{equation} \label{1}
c \frac{\partial}{\partial t} V(\bar{x},t) =
I(\bar{x},t) - V(\bar{x},t) + \int_{\Omega} K(|\bar{x}-\bar{y}|) S(V(\bar{y},t))d \bar{y},
\end{equation}
\[ t\in[0,T], \bar{x} \in \Omega \subset \Rb^2,\]
where the unknown $V(\bar{x},t)$ is a continuous function 
$V: \Omega \times [0,T] \rightarrow \Rb$, $I$, $K$ and $S$ are given functions; $c$
is a constant.
We search for a solution $V$ of this equation which satisfies
the initial condition
\begin{equation} \label{1a}
V(\bar{x},0)= V_0(\bar{x}), \qquad \bar{x}  \in \Omega.
\end{equation} 
Along with equation \eqref{1} we will also consider
\begin{equation} \label{2}
c \frac{\partial}{\partial t} V(\bar{x},t) =
I(\bar{x},t) - V(\bar{x},t)+ \int_{\Omega} K(|\bar{x}-\bar{y}|) S(V(y,t-\tau(\bar{x},\bar{y}))d\bar{y},
\end{equation}
\[  t\in[0,T],  \quad  \bar{x}  \in \Omega \subset \Rb^2,\]
where $\tau(\bar{x},\bar{y} )>0$ is a delay, depending on the spatial variables.
In the last case, the initial condition has the form
\begin{equation} \label{2a}
V(\bar{x},t)= V_0(\bar{x},t), \qquad \bar{x}  \in \Omega, \quad t \in [-\tau_{max},0],
\end{equation} 
where $\tau_{max}=\max_{\bar{x},\bar{y} \in \Omega} \tau (\bar{x},\bar{y})$.
By integrating both sides of \eqref{2} with respect to time on $[0,T]$, we obtain
the Volterra-Fredholm integral equation:
\begin{equation} \label{3}
c  V(\bar{x},t) = V_0(\bar{x})+ \int_0^t \left((I(\bar{x},s) - V(\bar{x},s)+ \int_{\Omega} K(|\bar{x}-\bar{y}|) S(V(y,s-\tau(\bar{x},\bar{y}))d\bar{y}\right) ds,  
\end{equation}
\[  t\in[0,T], \bar{x}  \in \Omega \subset \Rb^2.\]
The existence and uniqueness of solution of equation \eqref{1} was proved
in \cite{Pothast}, both in  the case of a smooth and discontinuous function $S$.
An analytical study of equation \eqref{2} was carried out in \cite{Faye},
where the authors have addressed the problems of existence, uniqueness and stability of solutions.

When solving numerically equations of the forms \eqref{1} and \eqref{2},
they are often reduced to the form \eqref{3}; therefore we begin by 
discussing literature on computational methods for Volterra-Fredholm equations.
 Starting with the one-dimensional case, without delay,
 Brunner has analysed the convergence of collocation methods \cite{Brunner}, while Kauthen
 has proposed continuous time collocation methods \cite{Kauthen}. In \cite{Han} an asymptotic error expansion for the Nystr\"{o}m method was proposed, which enabled the use
 of extrapolation algorithms to accelerate the convergence of the method.
 Another approach was developped by Z. Jackiewicz and co-authors \cite{Jacki1}, \cite{Jacki2}, who have applied Gaussian quadrature rules and interpolation
 the approximate the solution of integrodifferential equations
 modelling neural networks, which are similar to equation \eqref{1}.
 
 In the two-dimensional case, the required computational effort to solve
 equations  \eqref{1} and \eqref{3}  grows very fast as
 the discretization step is reduced, and therefore special attention has to be paid 
 to the creation of effective methods. 
 An important approach are the 
 low-rank methods, as those discussed in \cite{Xie}, when the kernel is approximated
 by polynomial interpolation,  which enables a significant reduction of the dimensions of the matrices. In \cite{Cardone}, the authors use an iterative method
 to solve linear systems of equations which takes into account the special form of the matrix to introduce parallel computation.
Concerning equation \eqref{2},  besides the existence and stability of
 solution, numerical approximations were obtained in \cite{Faye}.
The computational method applies a quadrature rule in space to reduce the problem
to a system of delay differential equations, which is then solved
by a standard algorithm for this kind of equations. A more efficient approach was recently
prpoposed in \cite{Hutt2010} \cite{Hutt2014}, where the authors introduce a new approach to deal with
the convolution kernel of the equation and use Fast Fourier Transforms to reduce 
signifficantly the computational effort required by numerical integration.

The above mentioned equations are known as Neural Field Equations (NFE) and  have played an important role in mathematical
neuroscience for a long time. Equation \eqref{1} was introduced first by Wilson
and Cowan \cite{WC}, and then by Amari \cite{Amari}, to describe excitatory
and inhibitory interactions in populations of neurons.
While in other mathematical models of neuronal interactions  the function
$V$ (membrane potential) depends only on time, in the case of NFE it is
a function of time and space.
The function $I$ represents external sources of excitation and
$S$ describes the dependence between the firing rate of the neurons and their membrane
potential. It can be either a smooth function (typically of sigmoidal type)
or a Heaviside function. The kernel function $K(|\bar{x}-\bar{y}|)$ gives the connectivity
between neurons in positions  $\bar{x}$ and $\bar{y}$. By writing the arguments of the function
in this form we mean that we consider the connectivity homogeneous, that
is, it depends only on the distance between neurons, and not on their specific
location.

According to many authors, realistic models of neural fields must take
into account that the propagation speed of neuronal interactions is finite,
which leads to NFE with delays of the form \eqref{2}. Here the delay $\tau$ is
space dependent, that is, it depends on the distance between the positions $x$ and $y$.

In \cite{LP1} a new numerical algorithm for the solution of equation \eqref{3} was introduced.
This algorithm combines  a second order implicit scheme for the space discretization 
with  Gaussian quadrature and collocation methods for the space discretization.
We refer to this article for the detailed description of this method and the analysis 
of its convergence.
In \cite{LP2} this algorithm was applied to the solution of 3 problems with physical interest.
The aim of the present paper is to present some more applications of the above mentioned algorithm
and to discuss its performance. 
\section{Numerical Simulations}
{\bf Example 1.} This algorithm is targeted directly to the application in Neuroscience and Robotics. As an illustration of such applications, we first consider a neural field described in \cite{Faye}, where the firing rate function has the form
\[
S(x)=\frac{2}{1+e^{- \mu x} },
\]
where $\mu \in \Rb^+$, and the connectivity function is given by
\[
K(r)= \frac{1}{\sqrt{2 \pi \xi_1^2}} \exp \left(-\frac{r^2}{2 \pi \xi_1^2} \right)-
\frac{A}{\sqrt{2 \pi \xi_2^2}} \exp \left(-\frac{r^2}{2 \pi \xi_2^2} \right),
\]
where $r=\|x-y\|_2=\sqrt{(x_1 -y_1)^2+(x_2 -y_2)^2}$ and $\xi_1, \xi_2, A  \in \Rb^+$.

It is   known \cite{Faye} that when there is no external input ($I(x,y,t)=0$), the stability of the trivial solution of this equation depends on the value of  $\mu ,$on the parameters of the connectivity  function and on the propagation speed. In particular, for each set of values of these parameters, it is possible to compute a bifurcation value $\mu_{bif}$, such that if  $\mu <\mu_{bif}$  the zero solution is  asymptotically stable, and otherwise it isn't.

We have considered some test cases, in order to analyse how the behaviour of solutions depends
on all these parameters. In each case, we have used our algorithm to approximate the solution with $h_t=0.5$ on the time interval $[0,T]$, where $T$ depends on the specific example (usually $25 \le T \le 75$). The parameters of the space discretization were $m=12$, $N=48$. We have used as initial condition $V_0(x) \equiv 0.01$
(in the delay case, we have  $V(t,x) \equiv 0.01$, if $t \in [-\tau_{max},0]$).
It is worth to note that according to \cite{Veltz} all the solutions of equation (\ref{2}) are bounded, so they can either tend to the trivial solution
($V\equiv 0$), to another steady state or just oscillate.

{\bf Case 1}. $\xi_1=0.1, \xi_2=0.2 , A=1$. In this case, the interaction between neurons is excitatory at short distances and inhibitory, at long ones
(so-called mexican hat connectivity).
For this case , we consider two different values of $\mu$: a) $\mu=15$, 
b) $\mu=45$; and two different forms of propagation: (i) no delay; ii) v=1.

The case $\mu=45$, without delay, was studied in \cite{Faye}. There the authors report that as time tends to infinity  the solution tends to a certain stationary state, with a biperiodical pattern. We have observed the same kind of behaviour, as can be seen in Fig.1. In this figure (as in other figures of the same type below) the graphs on the left-hand are the surface plots of the solution at equidistant moments in time ($t=30, 45,60,75$), while the graphs on the left-hand side are the corresponding contour plots. Note 
that the solution takes a long time to attain the stationary point. 
This can be observed in the graph of Fig.2. As other graphs of the same type below, in this graph the plots of the fucntions $V(x_1,t)$ and $V(x_2,t)$ are displayed, as functions of time, where $x_1$ and $x_2$ are  points located at the boundary and at the middle of the domain $\Omega$, respectively. Note that it is important that we apply an implicit method, which allows us to use a not very small stepsize in time. The fact that the solution does not tend to the zero stationary state indicates that this state is in this case unstable.
We have also studied the same case with propagation speed $v=1$. In this case  
the solution has a very different behaviour. It also tends 
to some non-zero stationary state, but obviously a different one (see Fig.3). 
When $\mu=15$, without delay,  the behaviour of the solution is typical to the case of a stable zero state: it tends rapidly  to 0 (see Fig. 5). Though during a certain time it exhibits a biperiodical pattern (see Fig.4), we observe that the amplitude decreases very fast with time. 
  
In the case $v=1$, the solution with $\mu=15$ has a similar behaviour to the no-delay case (see Fig.6). As could be expected, the main difference is that in this case the solution decays not so fast (see Fig.7).

\begin{figure}
\includegraphics[clip=true, trim= 1cm 6cm 1cm 6cm,  width=10cm]{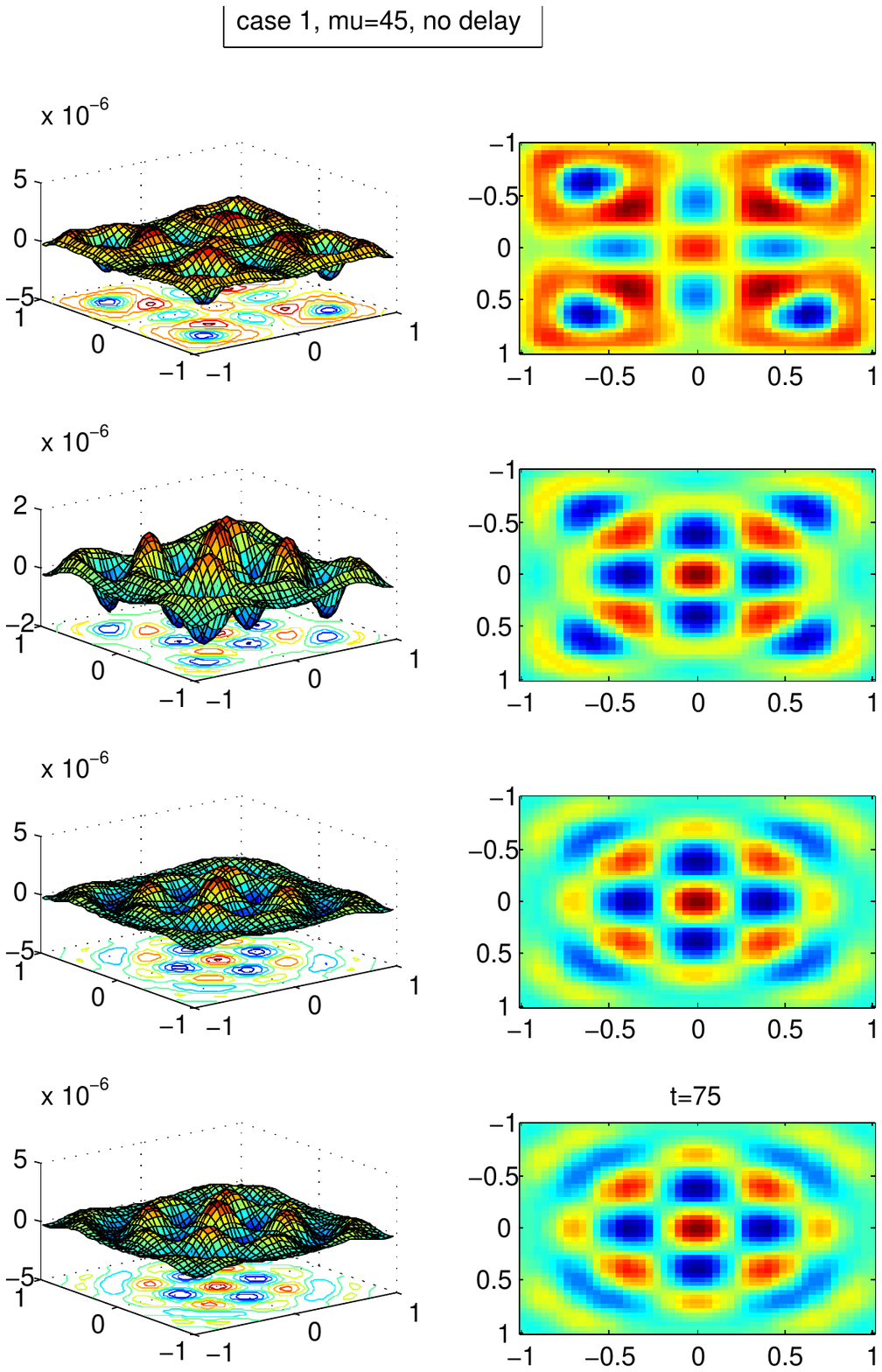}
\caption{ 3D-graphs (left) and contour plots (right) of the solution at $t=30,45,60,75$ , in the case 1, with $\mu=45$, with no delay.}
\end{figure}
\begin{figure}
\includegraphics[clip=true, trim= 1cm 9cm 1cm 9cm,  width=8cm]{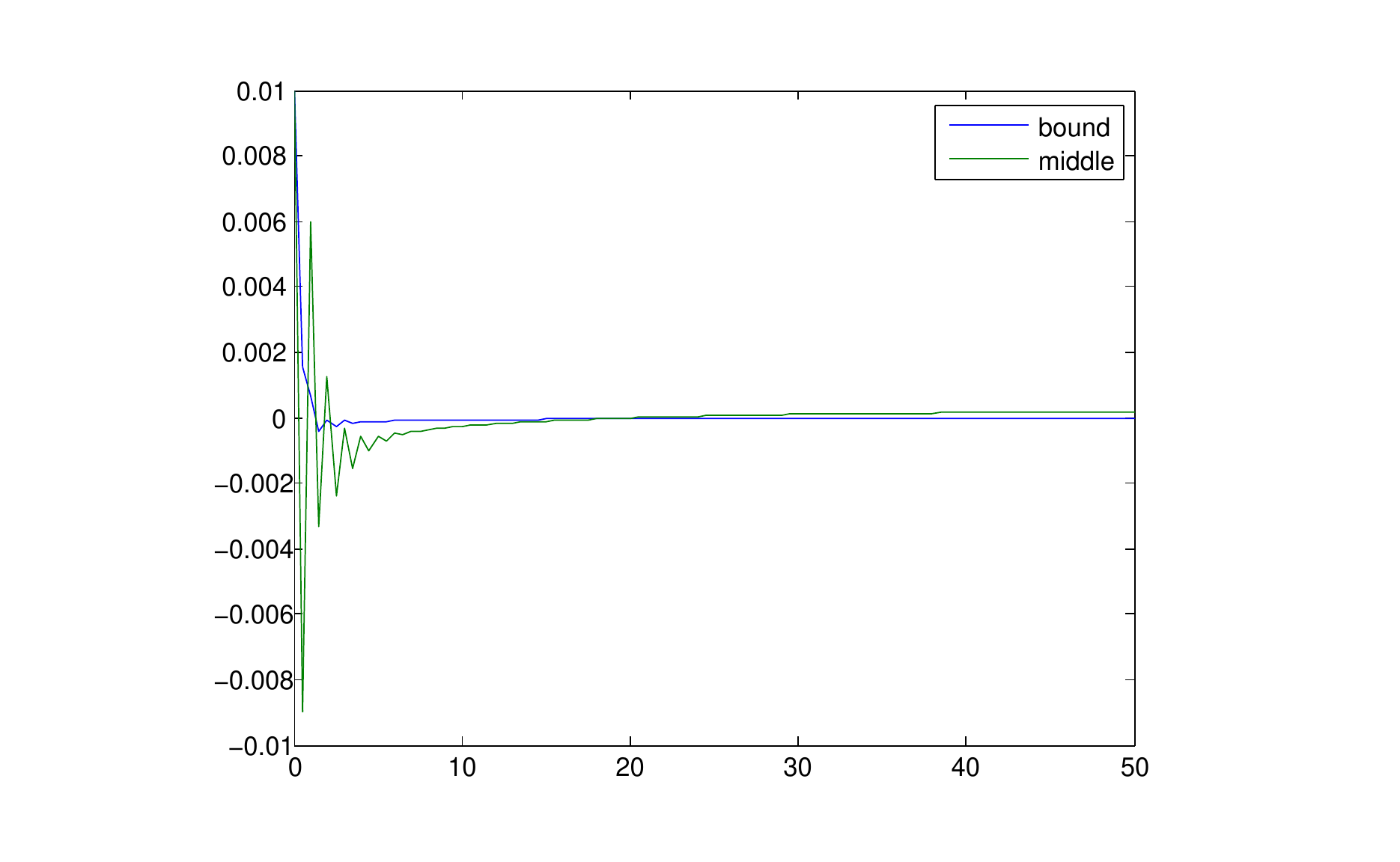} 
\caption{Graphs of $V(x_1,t)$  and $V(x_2,t)$ 
 in case 1,  $\mu=45$, with no delay}
\end{figure}
\begin{figure}
\includegraphics[clip=true, trim= 1cm 6cm 1cm 6cm,  width=10cm]{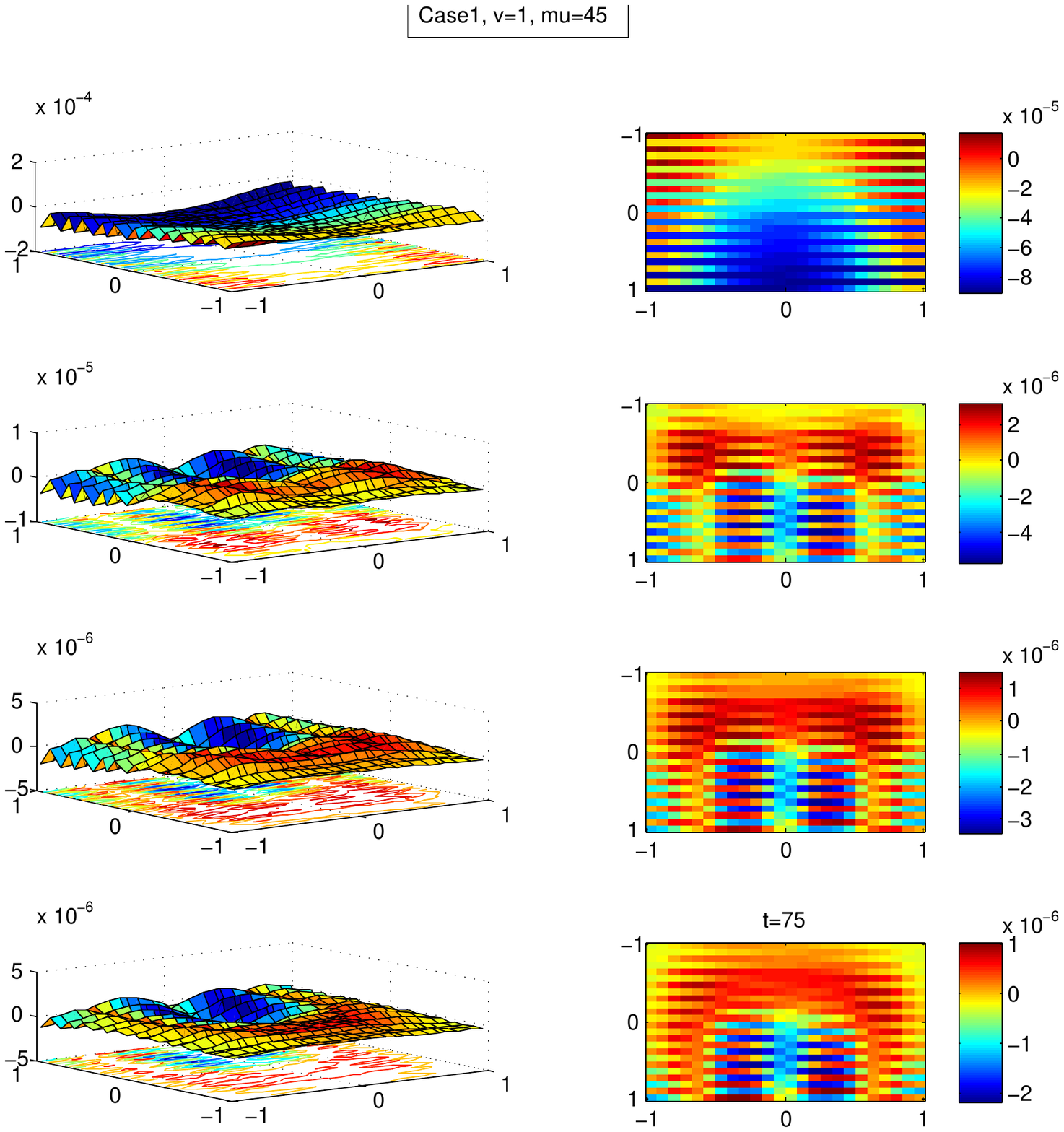}
\caption{ 3D-graphs (left) and contour plots (right) of the solution at $t=30,45,60,75$ , in the case 1, with $\mu=45$, with $v=1$.}
\end{figure}
\begin{figure}
\includegraphics[clip=true, trim= 1cm 6cm 1cm 6cm,  width=10cm]
{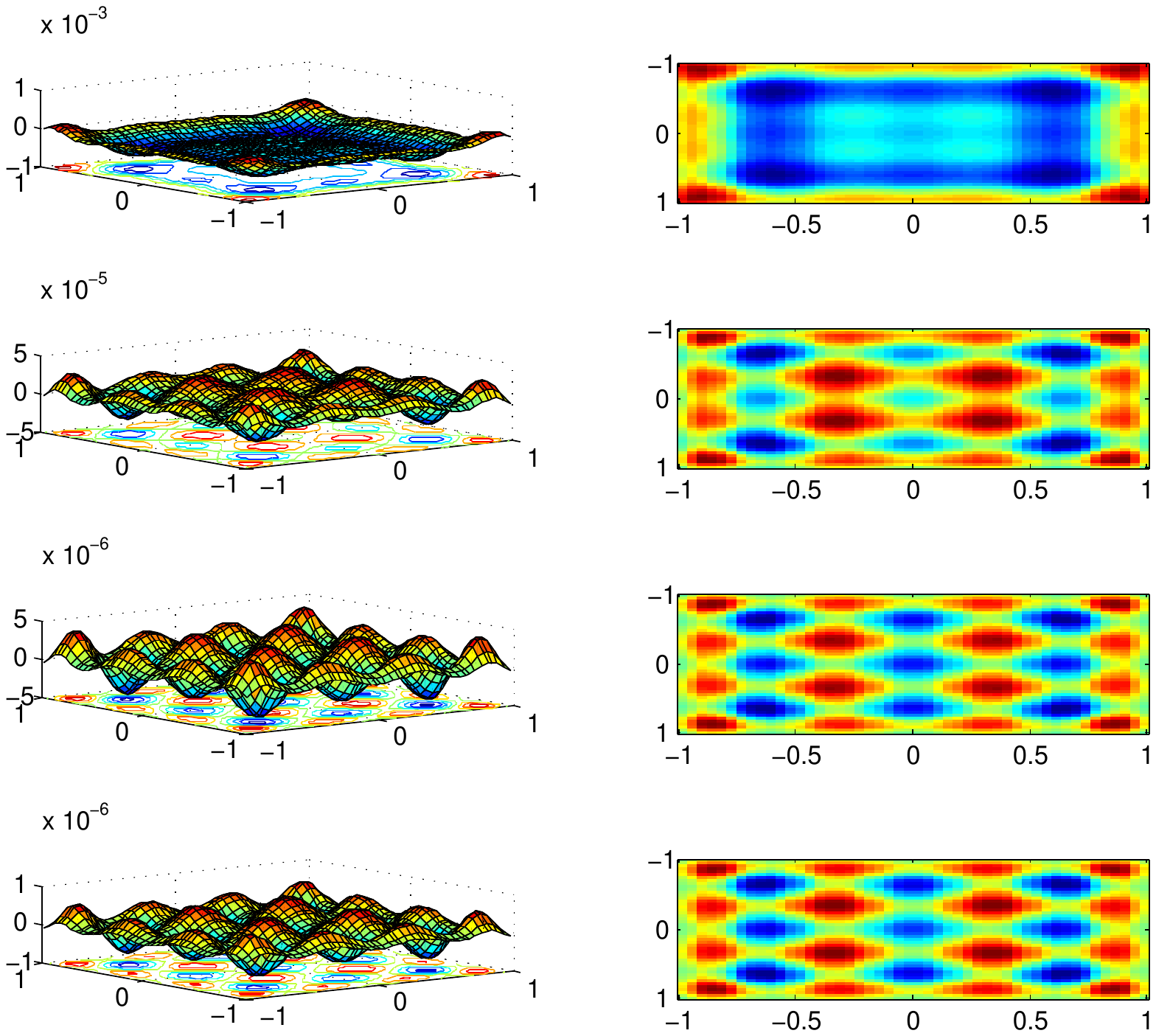}
\caption{ 3D-graphs (left) and contour plots (right) of the solution at $t=30,45,60,75$ , in the case 1, with $\mu=15$, with no delay.}
\end{figure}
\begin{figure}
\includegraphics[clip=true, trim= 1cm 9cm 1cm 9cm,  width=8cm]
{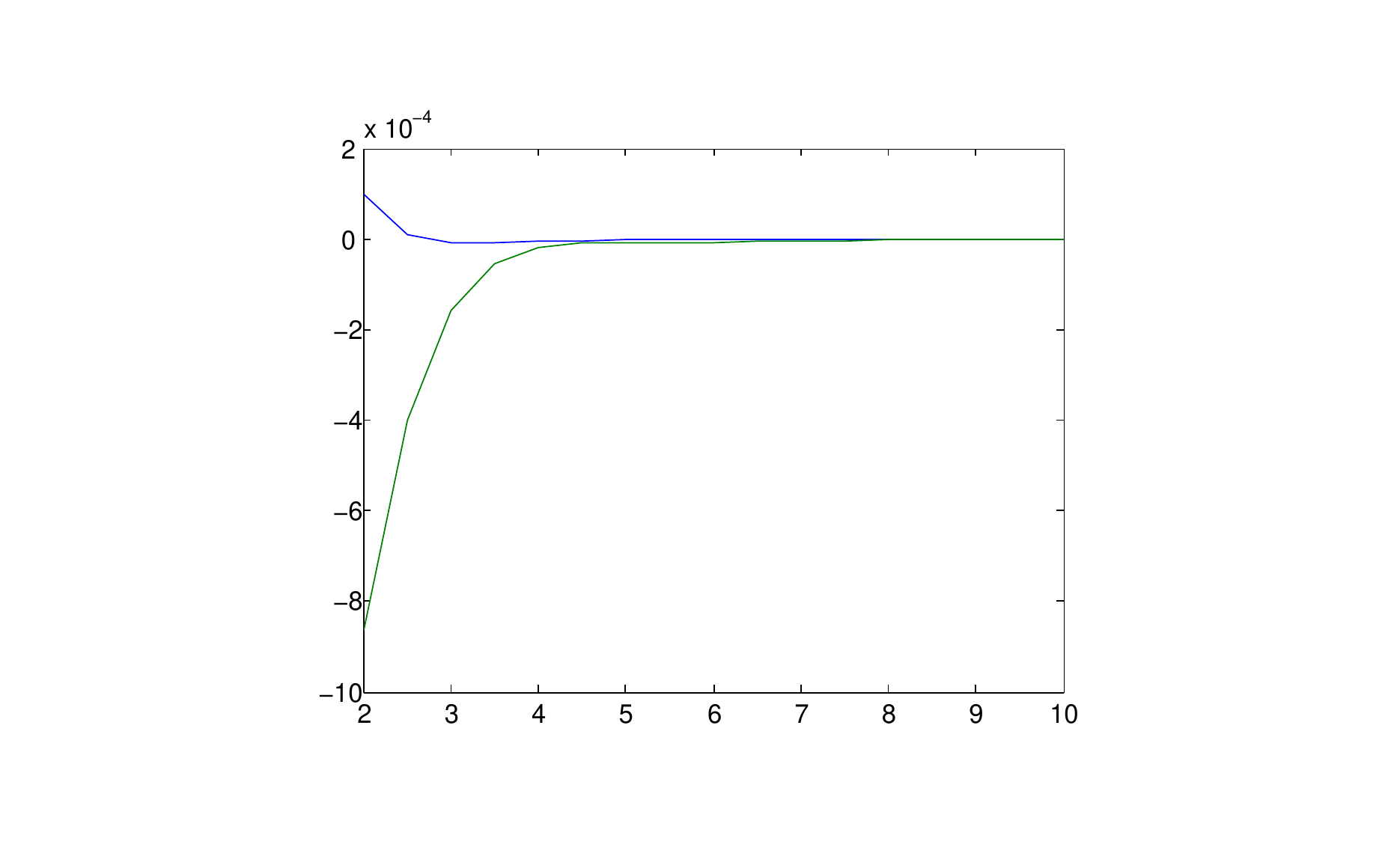} 
\caption{Graphs of $V(x_1,t)$  and $V(x_2,t)$ 
 in case 1,  $\mu=15$, with no delay.}
\end{figure}
\begin{figure}
\includegraphics[clip=true, trim= 1cm 6cm 1cm 6cm,  width=10cm]
{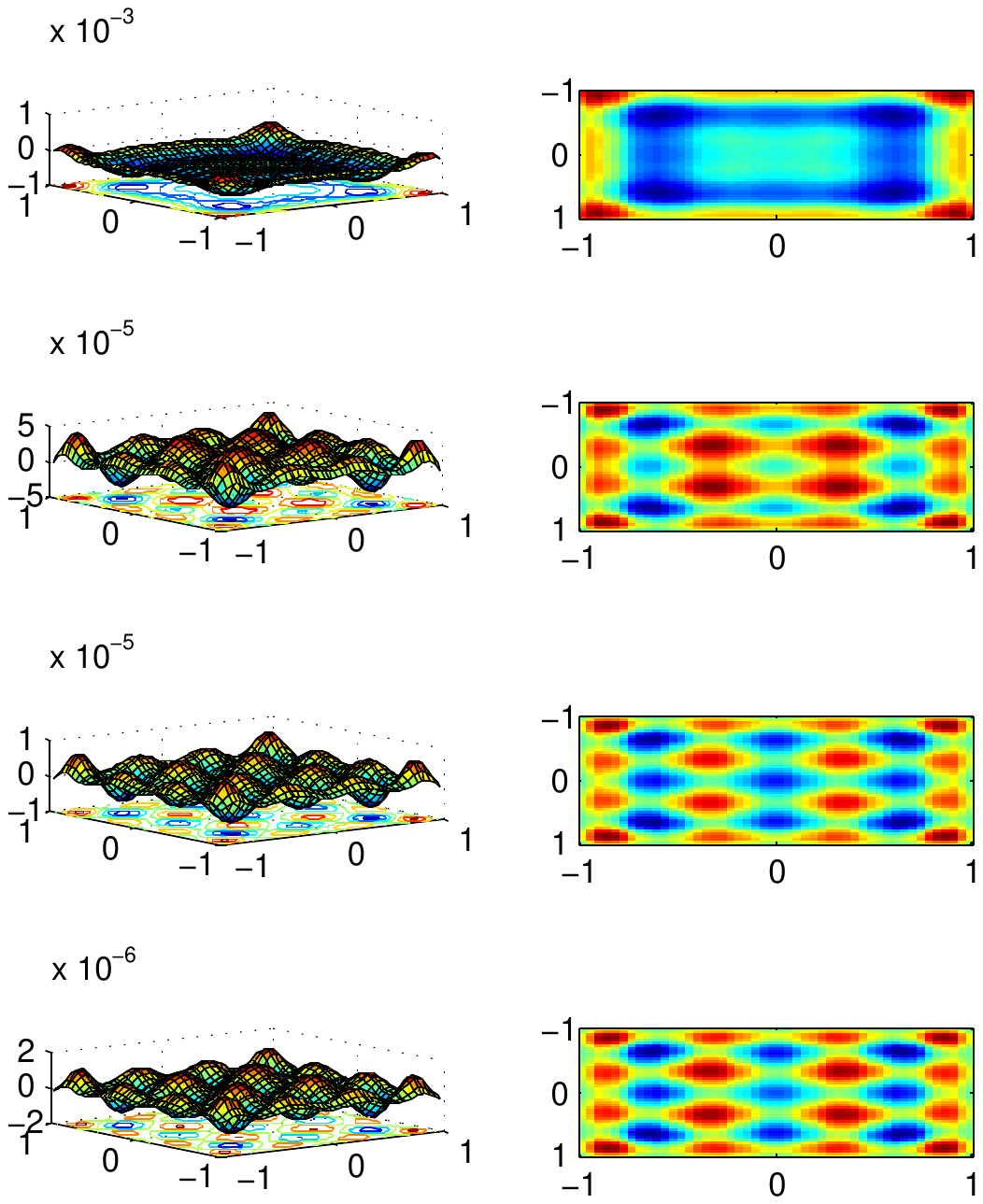}
\caption{ 3D-graphs (left) and contour plots (right) of the solution at $t=30,45,60,75$ , in the case 1, with $\mu=15$, with $v=1$.}
\end{figure}
\begin{figure}
\includegraphics[clip=true, trim= 1cm 9cm 1cm 9cm,  width=8cm]
{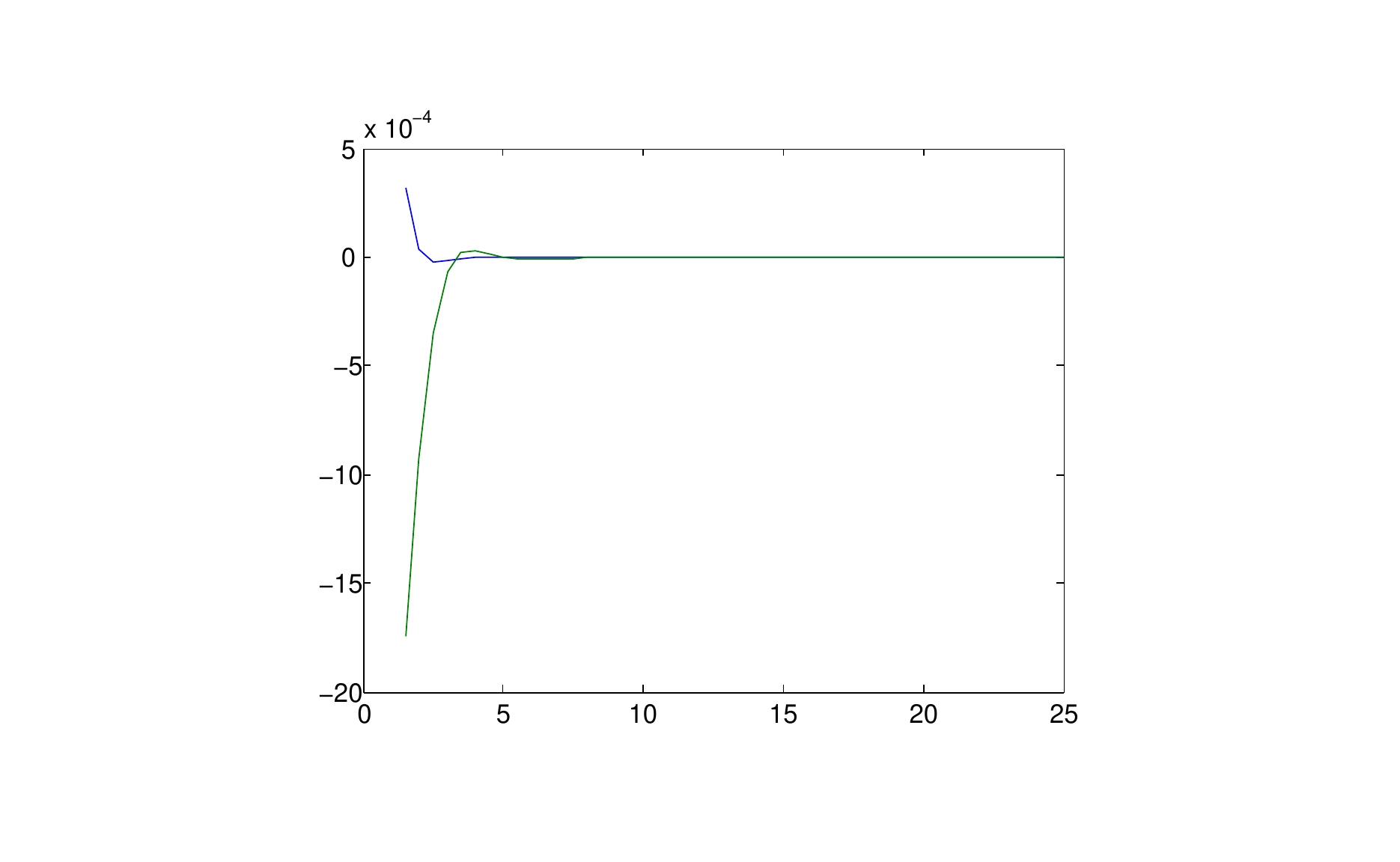} 
\caption{Graphs of $V(x_1,t)$  and $V(x_2,t)$ 
 in case 1,  $\mu=15$, with $v=1$.}
\end{figure}
{\bf Case 2}. Let $\xi_1=0.4, \xi_2=0.2 , A=1$ . In this case, we also have mexican-hat connectivity, but with a stronger excitatory component
For this example, we also consider two different values of  $\mu$: a)$\mu=10$, b) $\mu=15$; and two different forms of propagation: (i)no delay; ii)$v=0.1$.

When $\mu=15$, without delay,  the solution tends to a pattern, which is periodic with respect to $x$ and constant with respect to $y$ (see Fig. 8) . This suggests that the zero solution is unstable.
The evolution of the solution at a point in the boundary and another 
at the midlle is displayed in Fig. 9. In this graph we observe clearly that the solution tends to a nonzero stationary state.

The graphs of the solution for the case $\mu=15$, but with $v=0.1$, are
displayed in  Fig.10. We again observe convergence to a nonzero stationary state, which is confirmed by the graphs in Fig. 11. However, this is clearly not the same stationary state as in the non-delay case.

When $\mu=10$ we observe a quite different situation.
The graphs for this case,  with no delay  (see Fig.12), show clearly 
the convergence of the solution to the zero state. The evolution 
is plotted in the graphs of Fig.13.

Finally, in Fig. 14 the graphs of the solution are plotted in
the case $\mu=10$, $v=1$. The corresponding evolution graphs
are plotted in Fig. 15. In this case, the graphics also suggest that the solution tends to the zero state, but not so fast.

\begin{figure}
\includegraphics[clip=true, trim= 1cm 6cm 1cm 6cm,  width=10cm]
{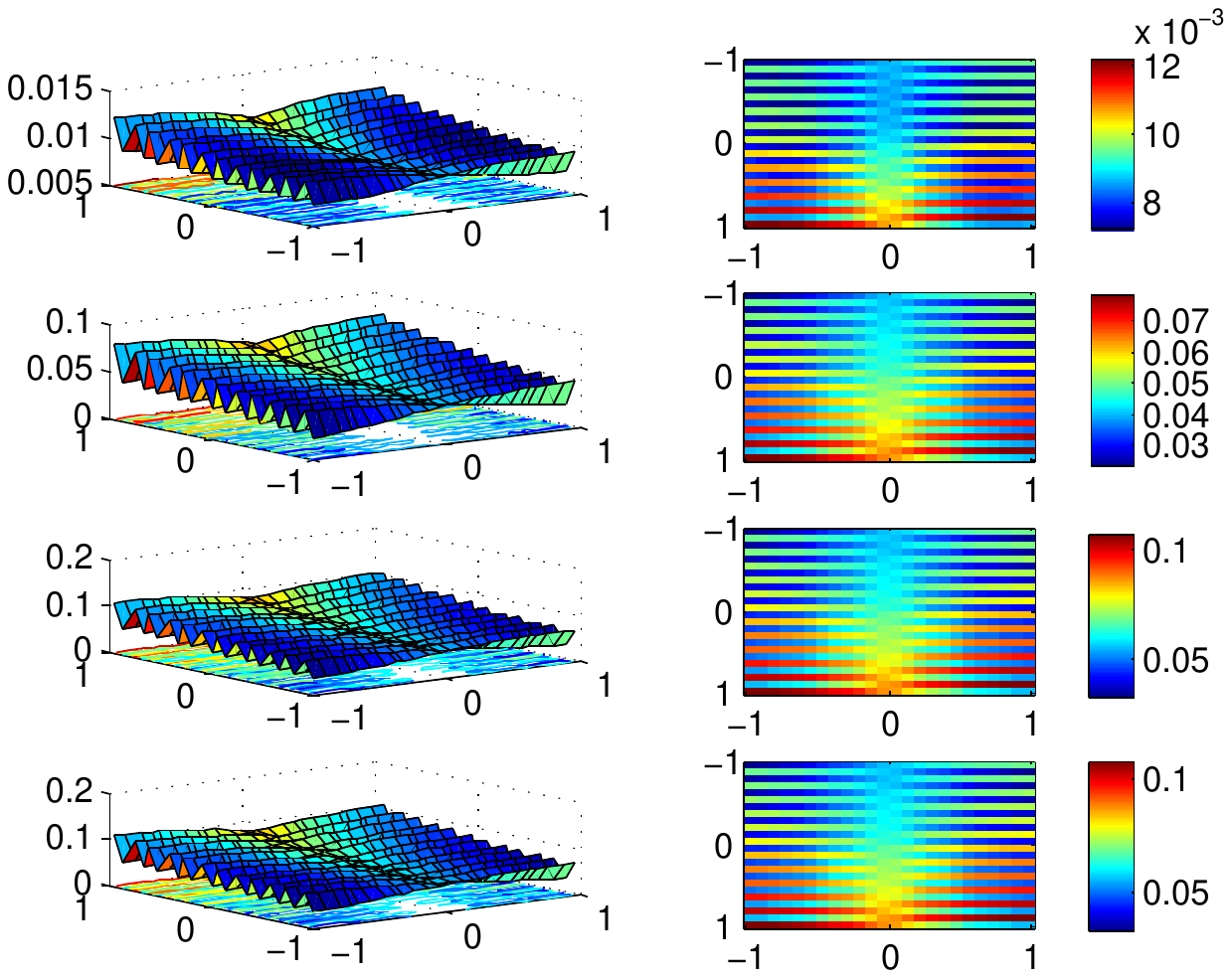}
\caption{ 3D-graphs (left) and contour plots (right) of the solution at $t=30,45,60,75$ , in the case 2, with $\mu=15$, with no delay.}
\end{figure}
\begin{figure}
\includegraphics[clip=true, trim= 1cm 9cm 1cm 9cm,  width=8cm]
{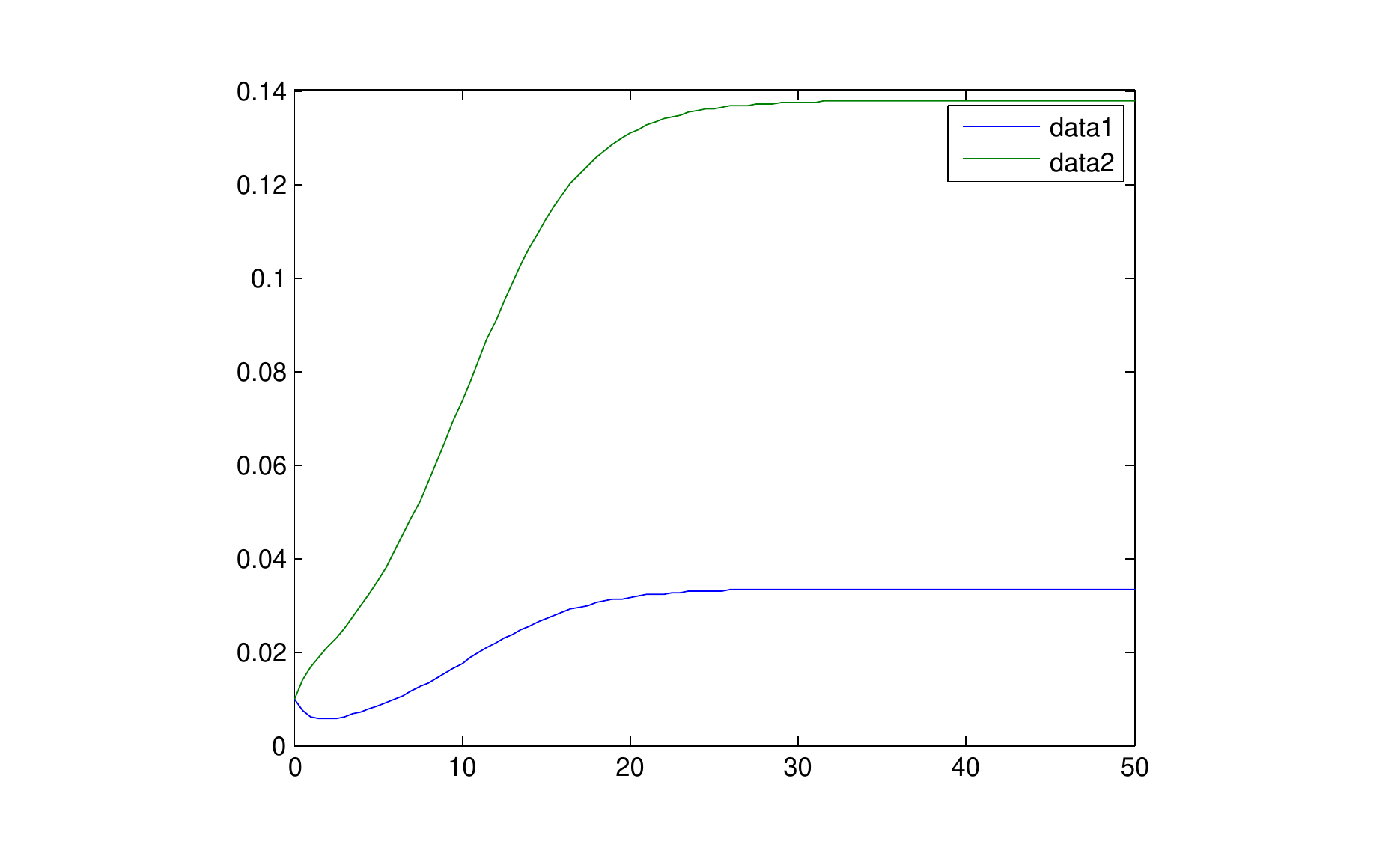} 
\caption{Graphs of $V(x_1,t)$  and $V(x_2,t)$ 
 in case 2,  $\mu=15$, with no delay}
\end{figure}
\begin{figure}
\includegraphics[clip=true, trim= 1cm 6cm 1cm 6cm,  width=10cm]{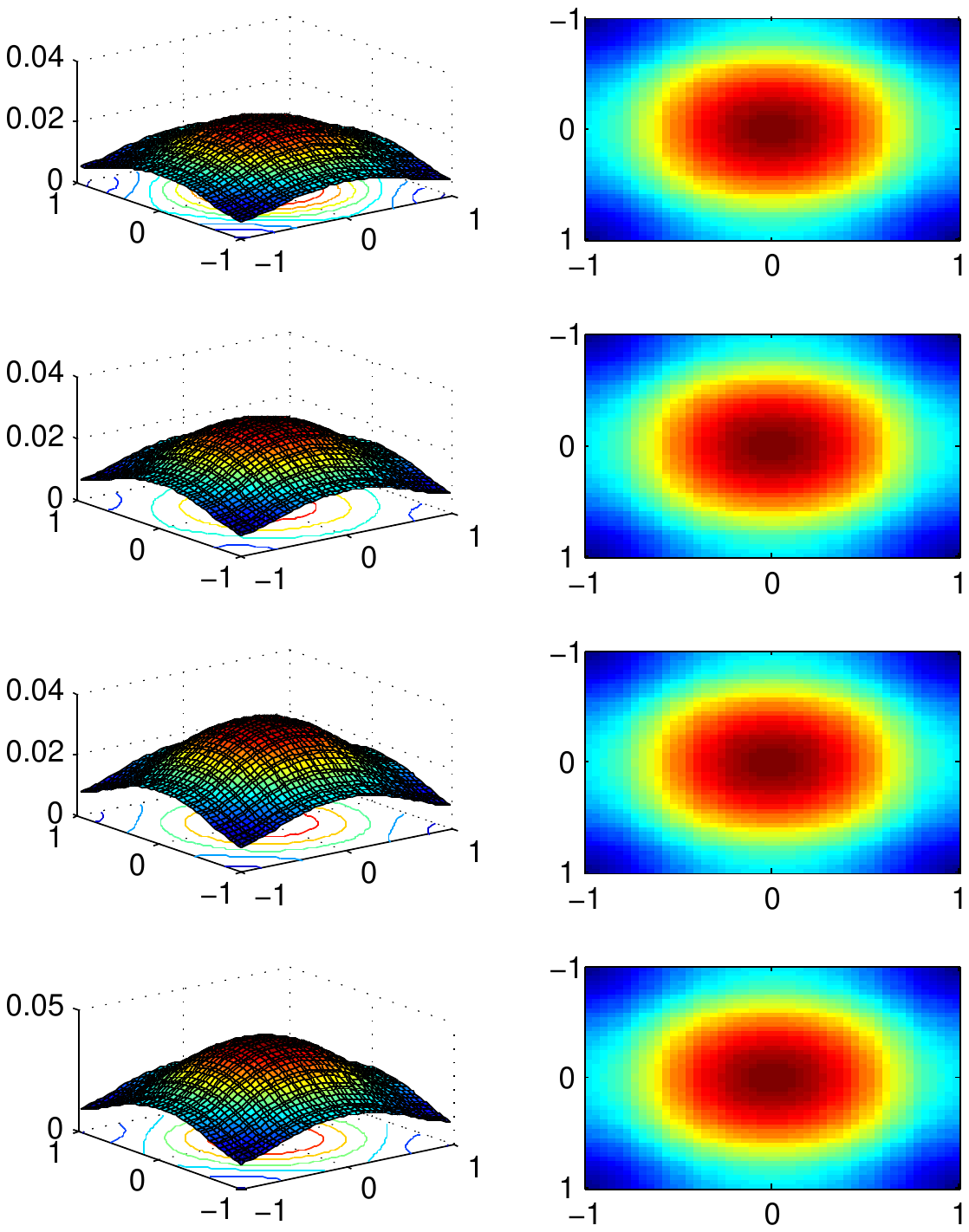}
\caption{ 3D-graphs (left) and contour plots (right) of the solution at $t=30,45,60,75$ , in the case 2, with $\mu=15$, with $v=0.1$.}
\end{figure}
\begin{figure}
\includegraphics[clip=true, trim= 1cm 9cm 1cm 9cm,  width=8cm]
{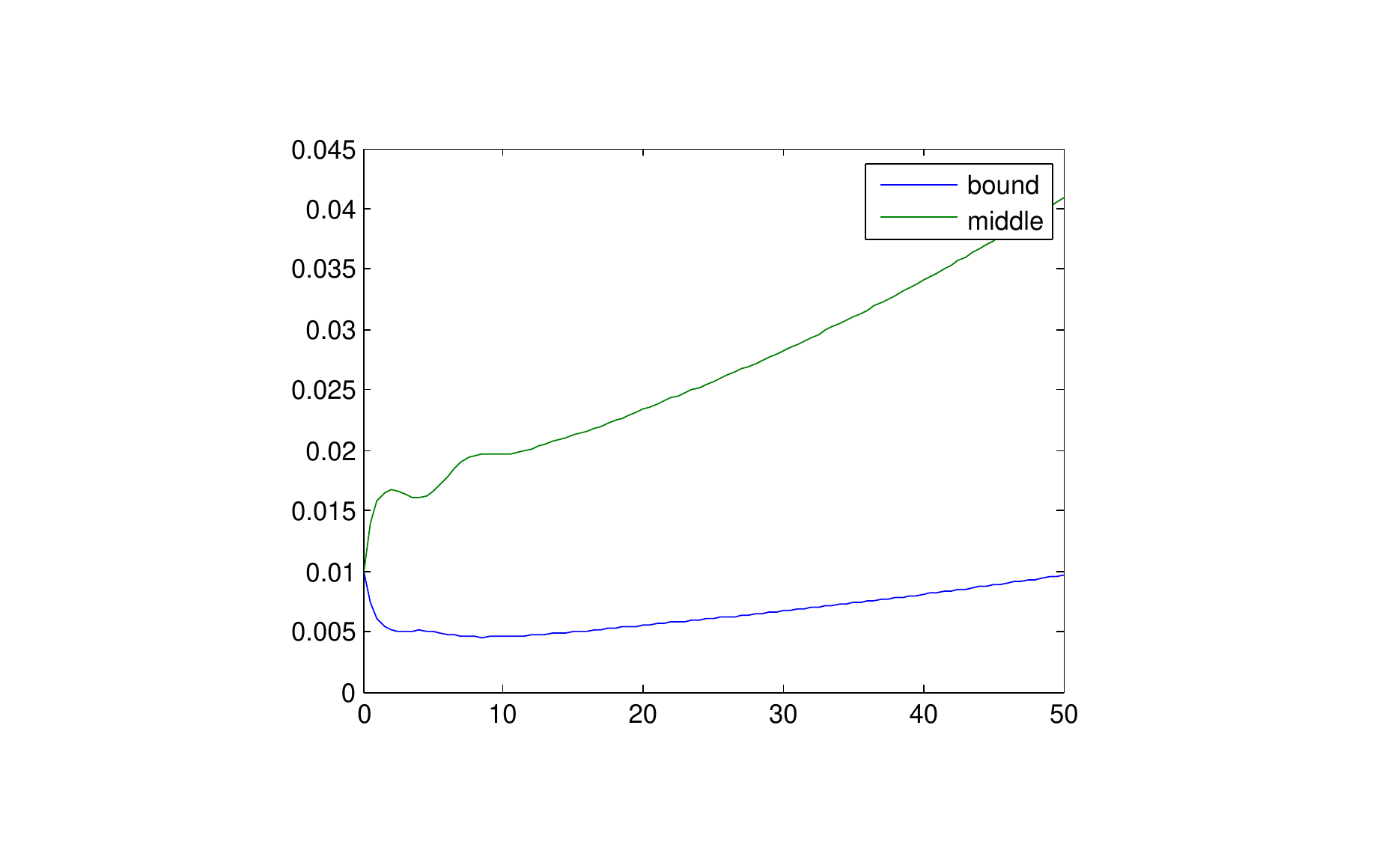} 
\caption{Graphs of $V(x_1,t)$  and $V(x_2,t)$ 
 in case 2,  $\mu=15$, with $v=0.1$.}
\end{figure}
\begin{figure}
\includegraphics[clip=true, trim= 1cm 6cm 1cm 6cm,  width=12cm]
{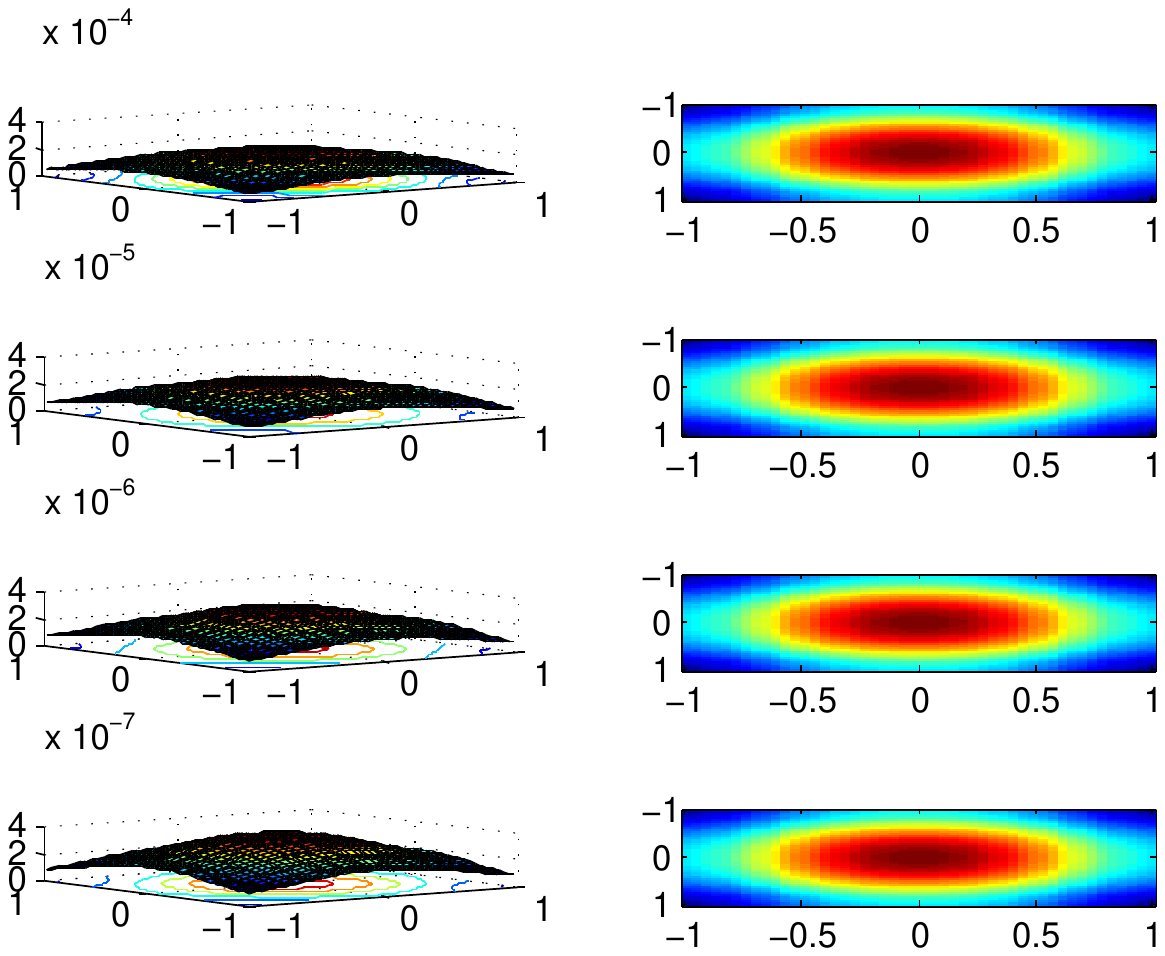}
\caption{ 3D-graphs (left) and contour plots (right) of the solution at $t=30,45,60,75$ , in the case 2, with $\mu=10$, with no delay.}
\end{figure}
\begin{figure}
\includegraphics[clip=true, trim= 1cm 9cm 1cm 9cm,  width=8cm]
{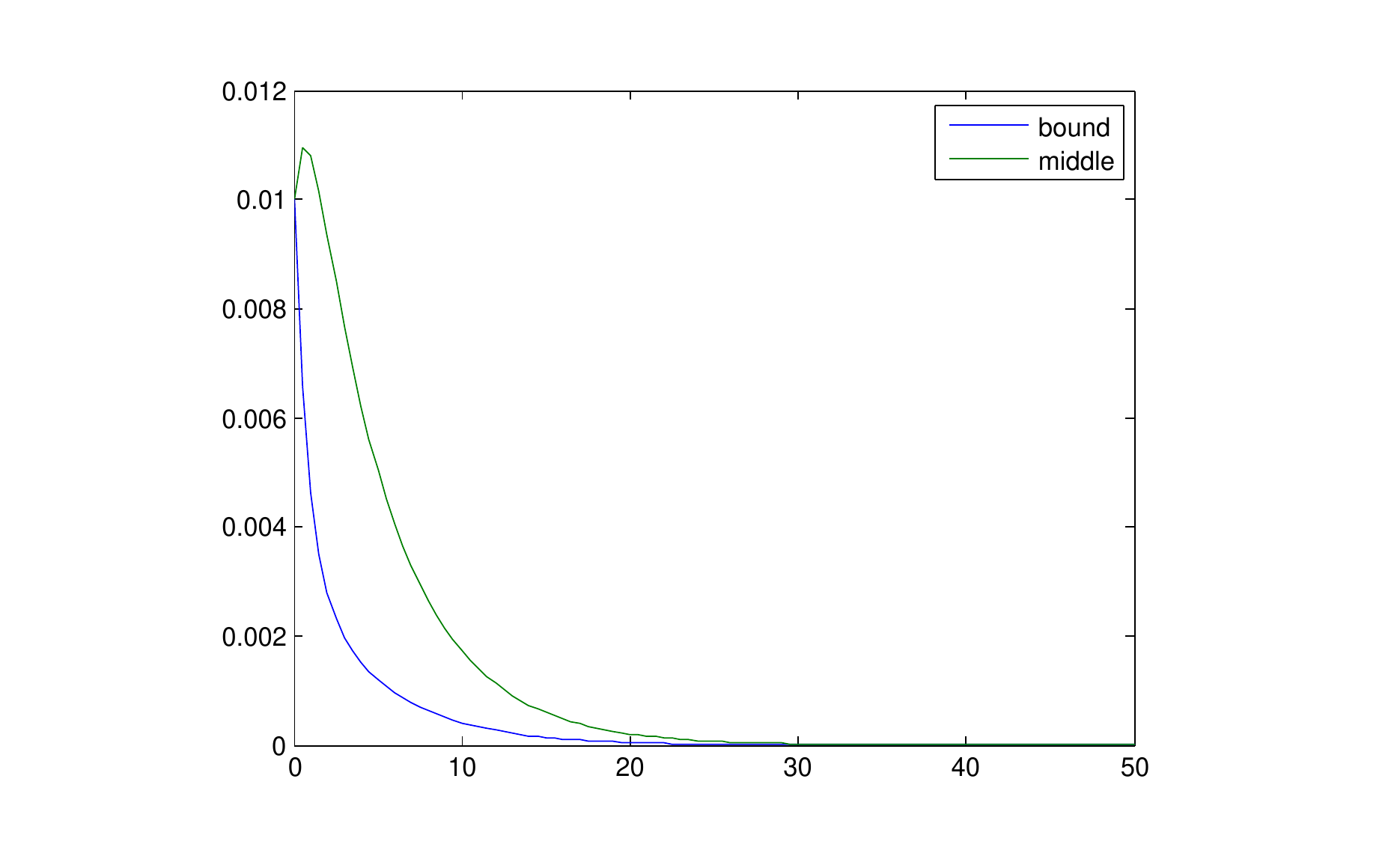} 
\caption{Graphs of $V(x_1,t)$  and $V(x_2,t)$ 
 in case 2,  $\mu=10$, with no delay.}
\end{figure}
\begin{figure}
\includegraphics[clip=true, trim= 1cm 6cm 1cm 6cm,  width=12cm]
{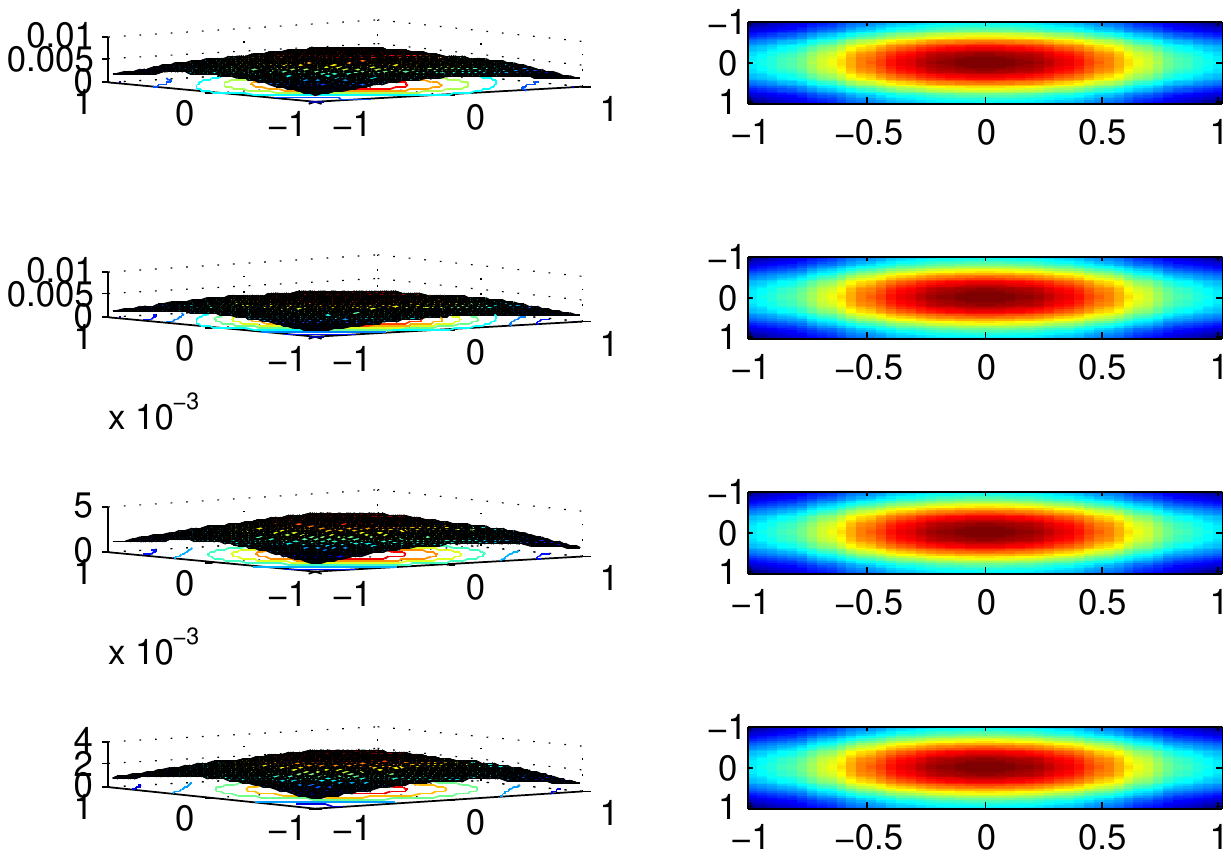}
\caption{ 3D-graphs (left) and contour plots (right) of the solution at $t=30,45,60,75$ , in the case 2, with $\mu=10$, with $v=0.1$.}
\end{figure}
\begin{figure}
\includegraphics[clip=true, trim= 1cm 9cm 1cm 9cm,  width=8cm]
{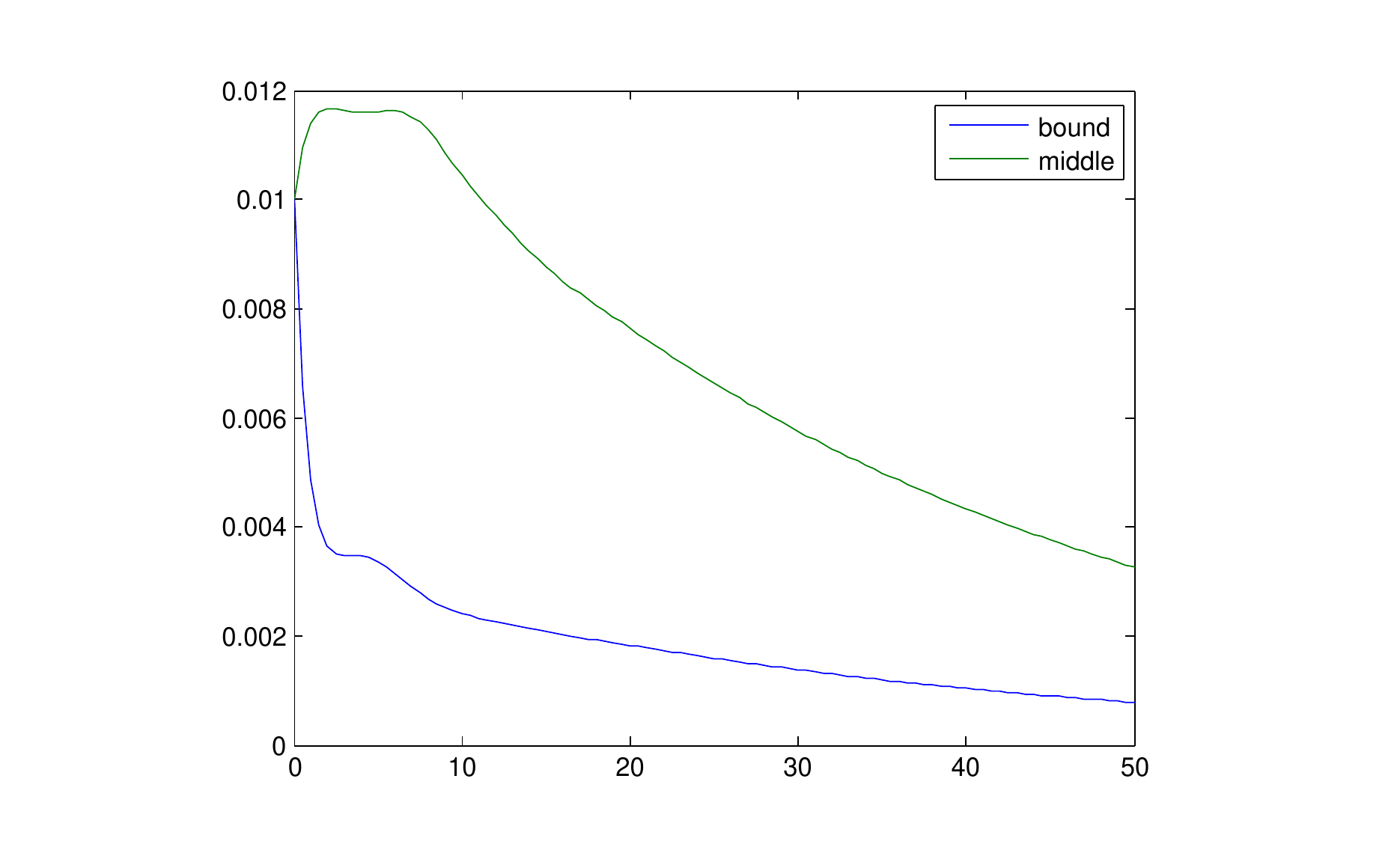} 
\caption{Graphs of $V(x_1,t)$  and $V(x_2,t)$ 
 in case 2,  $\mu=10$, with $v=0.1$.}
\end{figure}
\begin{table}
\begin{displaymath}
\begin{array} {|c|c|c|c|c|}
\hline
\textstyle{Case} &  \mu & v & \textstyle{Type \, of \, Solution}  & \textstyle{Figures} \\
\hline
1  &  45 & \textstyle{no \,delay} &  \textstyle{Convergence \, to  \,a  \, nontrivial \, state } & 1,2  \\
1 &  45 &  1  &  \textstyle{Convergence \,to \, a \, nontrivial state } &  3 \\
1 &  15 &  \textstyle{no \, delay}  &  \textstyle{Convergence to the zero solution } &  4,5 \\
1 &  15 &  1  &  \textstyle{Convergence \, to \,the \,zero \, solution} &  6,7 \\
2 &  15 &  \textstyle{no \,delay}  &  \textstyle{Convergence \,to \, a \,nontrivial \,state } &  8,9\\
2 &  15 &  0.1  &  \textstyle{Convergence \, to \,a \,nontrivial \, state} &  10,11 \\
2 &  10 &  \textstyle{no\, delay}  &  \textstyle{Convergence \, to \, the \, zero \, solution} &  12,13 \\
2 &  10 &  0.1  &  \textstyle{Convergence \,to \, the \, zero \,solution} &  14,15 \\
\hline
\end{array}
\end{displaymath}
\caption{ Summary of the numerical results for Example~1. }
\end{table}
The results of the numerical simulations concerning Example~1 are summarized in Table~1. 
Both in Case 1 and Case 2, when the value of $\mu$  is higher than a certain threshold we observe the 
convergence of the solution to  certain non-trivial stationary states, while for lower values of $\mu$ 
the solution converges to the zero state, which in this case is stable. This is in agreement with the
results presented  in \cite{Faye}, where it is mentioned that there is exists a bifurcation value for $\mu$.
Concerning the influence of propagation speed, we observe that in the case of a stable zero state
the solution converges to this state also in the case of a finite propagation speed, though the convergence
 is slower in this case (as it could be expected). In the case of the nontrivial stationary states,
the  solution converges to different states, depending on the propagation speed.
The main difference between Case~1 and Case~2 is that in Case 2 the critical value of $\mu$ (for which the
zer state becomes unstable) is lower. We recall that both in Case~1 and Case~2 we have mexican hat 
connectivity, but in Case~2 the parameters $\xi_1$ and $\xi_2$ are higher.

	{\bf Example 2.} We are now concerned about a neural field, which is described in \cite{Hutt2014}. In this case we
	have a connectivity kernel of the form: 
	\begin{equation} \label{ker2}
  K( \bar{x}, \bar{y})= K_0 \sum_{i=0}^2 \cos (\bar{k}_i . \bar{x}-\bar{y}) 
	\exp \left(- \frac{\|\bar{x}-\bar{y}\|}{\sigma} \right),
	\end{equation}
where $\bar{k}_i=k_c(\cos(\phi_i), \sin(\phi_i)) $, $\phi_i=i \pi/3$ and $K_0$, $k_c$ and $\sigma$ are
positive constants. In our numerical simulations we have used $K_0=1.5$, $\sigma=10$, $k_c=\pi$. This kernel reflects  spatial hexagonal connections, which have been found in the visual cortex of monkeys
\cite{Lund} (see fig. 16).  
In order to test our program, we have reproduced the numerical simulation carried out by Hutt and Rougier in \cite{Hutt2014}.
The firing rate function in this case has the form of the sigmoidal function
\[
S(x)= \frac{2}{1+\exp(-5.5(x-3))}
\]
and the external input is given by
\[
I(\bar{x},t)= I_0+ \frac{1}{ \pi \sigma_I^2} \exp (-\bar{x}^2 / \sigma_I^2),
\]
with $I_0=2$ and $\sigma_I=0.2$. Note that this external input is localized at the domain center, so that the function $I$ has very high derivatives 
in the neighborhood of this point. This property is reflected by the solution, which (as can be seen in fig. 17) has 
a very thin peek at the  domain center.
This creates some difficulties to the application of
our program. Namely, we cannot apply the rank reduction technique, described in \cite{LP1}, to reduce
the size of the matrices. It is not possible in this case to approximate accurately the solution by an interpolating polynomial
with a small number of interpolation points. We have overcome this difficulty by applying our algorithm without rank reduction.
Actually, since the discretization method converges fast, it is possible to obtain an accurate approximation of the solution
with a moderate number of gridpoints and a large stepsize in time (in our simulation, we have used $N=48$
and $ht=0.08$). In this way, it is possible
to carry out numerical simulations with a not too high computational effort.
The resolution of the solution graph can then be improved, obtaining additional solution values by interpolation.
With this purpose, we have used cubic interpolating polynomials at each interval $[x_i, x_{i+1}]$.

Just as reported in \cite{Hutt2014}, we note that the activity of the neural field is first concentrated in the domain center and then propagates to the boundary, with a sequence of peeks, which exhibit the same hexagonal pattern as the 
connectivity function. 

The aim of the numerical experiments concerning this example is to investigate the influence of the 
delay $\tau(x,y)$ on the solution of equation (\ref{2}). As described in Sec. 1, we assume that this delay 
results from the time spent by electric signals to travel between neurons. Therefore, it has the form
\[ \tau(x,y)= \|x-y\|/v,
\]
where $v$ is the velocity of signal propagation (which is supposed to be constant in time and uniform
in space).

Figures 17-20 illustrate the case where $v=10$.  In each figure, we have a 3D graphic of the solution 
(left-hand side) and the corresponding contour plot 
(right-hand side). As we can see from these figures, at each instant  of time,
it is possible to identify a region of activity, that is, a region $\Omega_a(t)$, such that if
$(x,y) \in  \Omega_a(t)$, then $V(x,y,t) \ge V_{th}$, where $V_{th}$ is some threshold value of the potential , under which we consider that there is no significant activity. In
our experiments, we have considered $V_{th}=2.1$, which was found experimentally and seems to be a good criterium.. As mentioned before, this domain   has the approximate
form of a hexagon (reflecting the characteristics of the connectivity function (\ref{ker2})).
We define the radius of activity $r_a(t)$ by
\[ r_a(t)=\max_{(x,y) \in \Omega_a(t)} \|(x,y)\|_2 .
\] 
The figures show clearly that $r_a(t)$ is an increasing function, and our purpose is to use the numerical results to analyse this function and investigate how it depends on the delay and, in particular, on the velocity $v$.

We note also that in the case of {\bf no delay} (which corresponds to infinite velocity $v$) the solution of the problem  in the initial period ($t \le 0.24$) has much higher values in the center than 
in the periphery of the domain $\omega$ (as it happens in the delay case). However in this case there is no sharp transition
between activity and non-activity  zones, and therefore we cannot identify clearly a domain of activity   $\Omega_a(t)$.
This case is illustrated by figure 21, where the contour plot of the solution is displayed in
the cases of $t=0.16$ and $t=0.32$.

In Table~2, we display the values of $r_a(t)$ for $0<t \le 0.96$ in three different cases: $v=10$, $v=20$ and no delay.
Some interesting conclusions can be deduced from this table.
\begin{itemize}
\item In the first moments after the start of stimulation ($t \le 0.16$), the velocity $v$ seems not to influence
the activity radius $r_a(t)$ . This is not surprising, because in the first moments the activity of the field  is significant 
only in a small region and the effect of delay within this region very is small (the time spent by signals to travel
within this region does not play an important role).  The radius of activity in this case depends mainly on the properties of 
the connectivity kernel, which does not depend on $v$.
\item  However, as time passes, the situation changes. For $t=0.24$, we see that
in the case of no delay the radius of activity is $r_a=14.142= l \sqrt{2}$, which means that from this moment there is activity in the whole domain.   In the delay case, $r_a (t)$ continues to grow with $t$, which physically means
that there is a wave front moving from the center to the periphery of the domain.  
\item This wave front
does not move with a constant speed (the differences between $r_a(t_2) -r_a(t_1)$ are not proportional to the
 time interval lengths ($t_2-t_1)$.
\item However, when we compare the cases $v=10$ and $v=20$ (for $t \ge 0.4$) , we see that the average propagation
speed of the wave front (defined as $r_a(t)/t$ ) is nearly the double in the second case, compared with the first one.
This means that for sufficiently large time, the transmission delay is the main factor for the definition of the activity radius.
Again that can be explained by the fact that when we have a large activity domain, the  transmission time between the neurons
has a strong influence in the dynamics of the neural field.   
\end{itemize}

\begin{figure}[p]
\centering
\includegraphics[clip=true, trim= 2cm 10cm 2cm 10cm,  width=9cm]{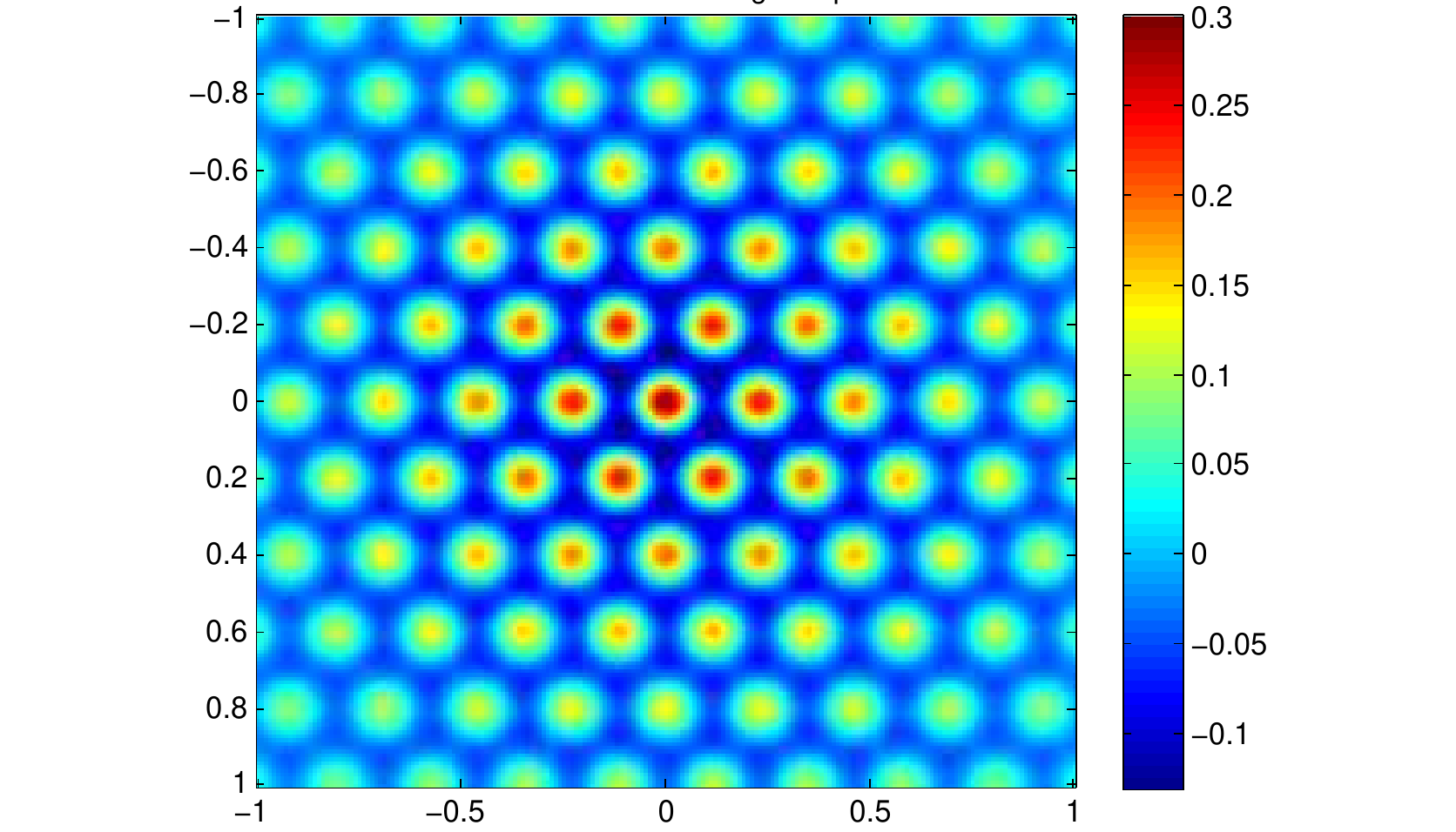} 
\caption{Plot of the connectivity function $K(x,y)$  of example~2, see formula (\ref{ker2})}.
\end{figure}

The numerical results were obtained from a numerical simulation using the method under discussion, with $v=10$, with stepsize $h_t=0.08$ and $N=48$.
In each figure, the contour plot displayed on the right-hand side corresponds to the 3D-graphic given on the left-hand side.
\begin{figure}[p]
\begin{displaymath}
\begin{array}{cc}
\includegraphics[clip=true, trim= 2cm 9cm 2cm 9cm,  width=8cm]{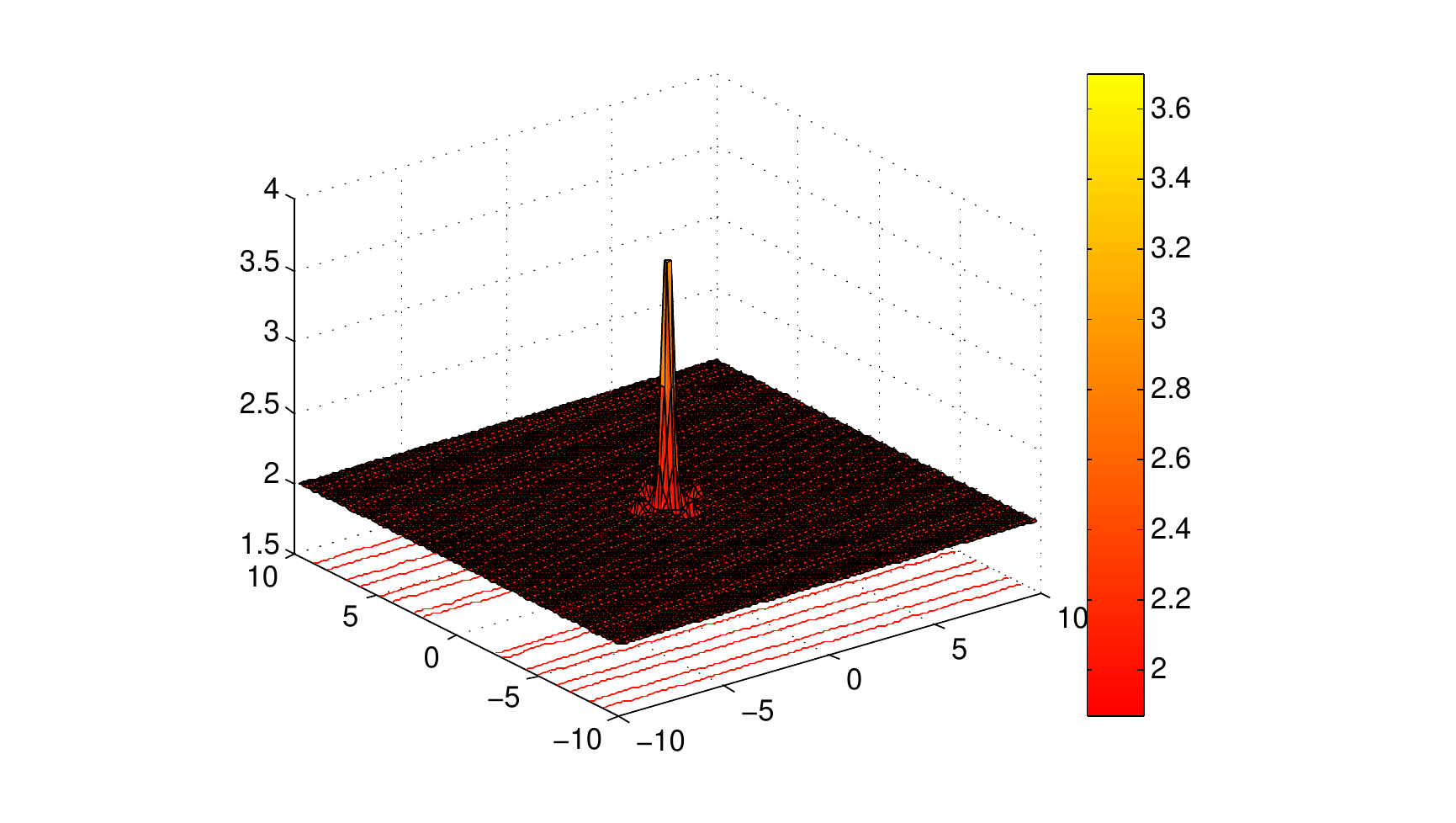} & 
\includegraphics[clip=true, trim= 2cm 9cm 2cm 9cm,  width=8cm]{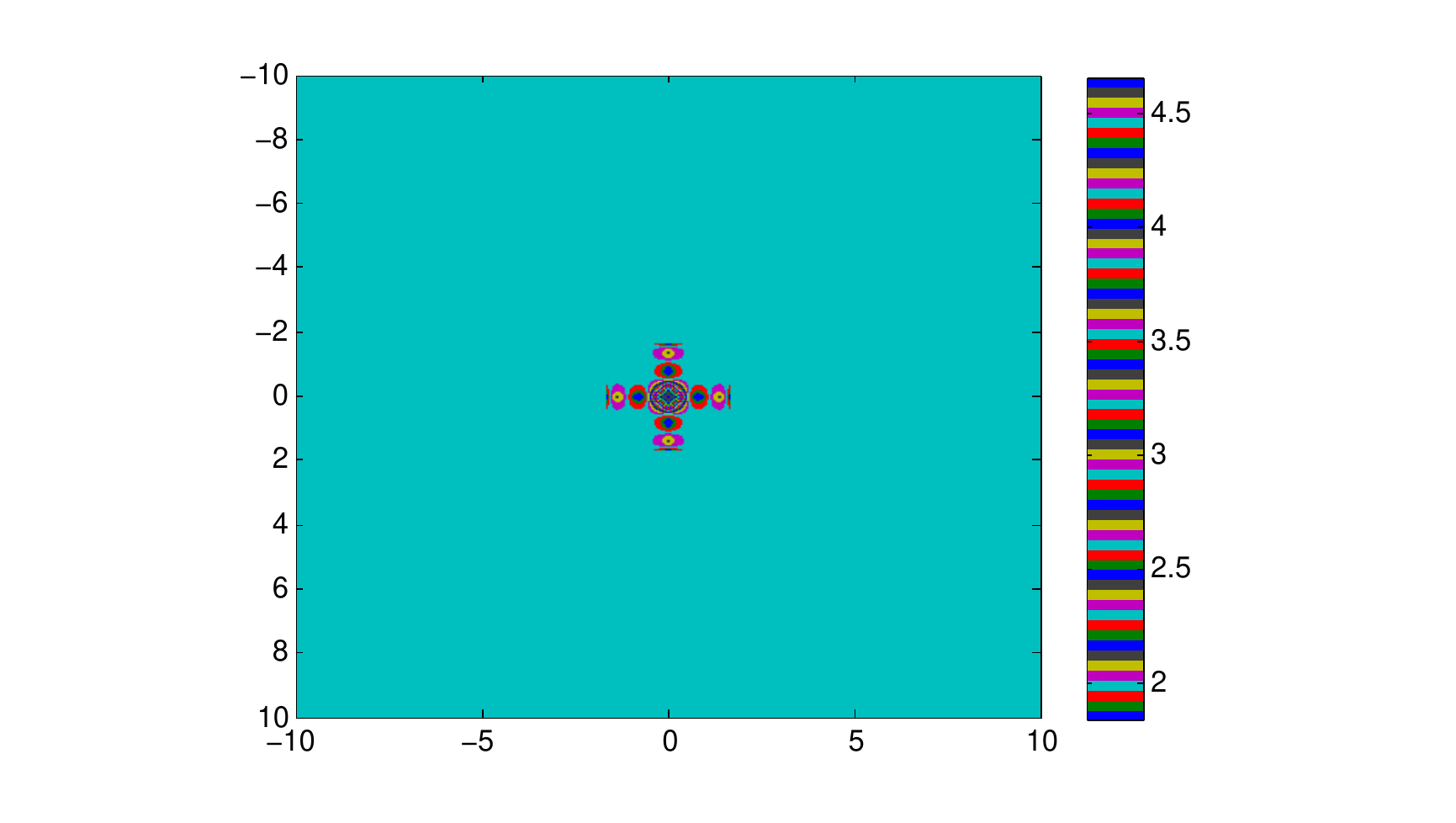} 
\end{array}
\end{displaymath}
\caption{Solution at $t=0.16$, with $v=10$.}
\end{figure}
\begin{figure}[p]
\begin{displaymath}
\begin{array}{cc}
\includegraphics[clip=true, trim= 2cm 9cm 2cm 9cm,  width=8cm]{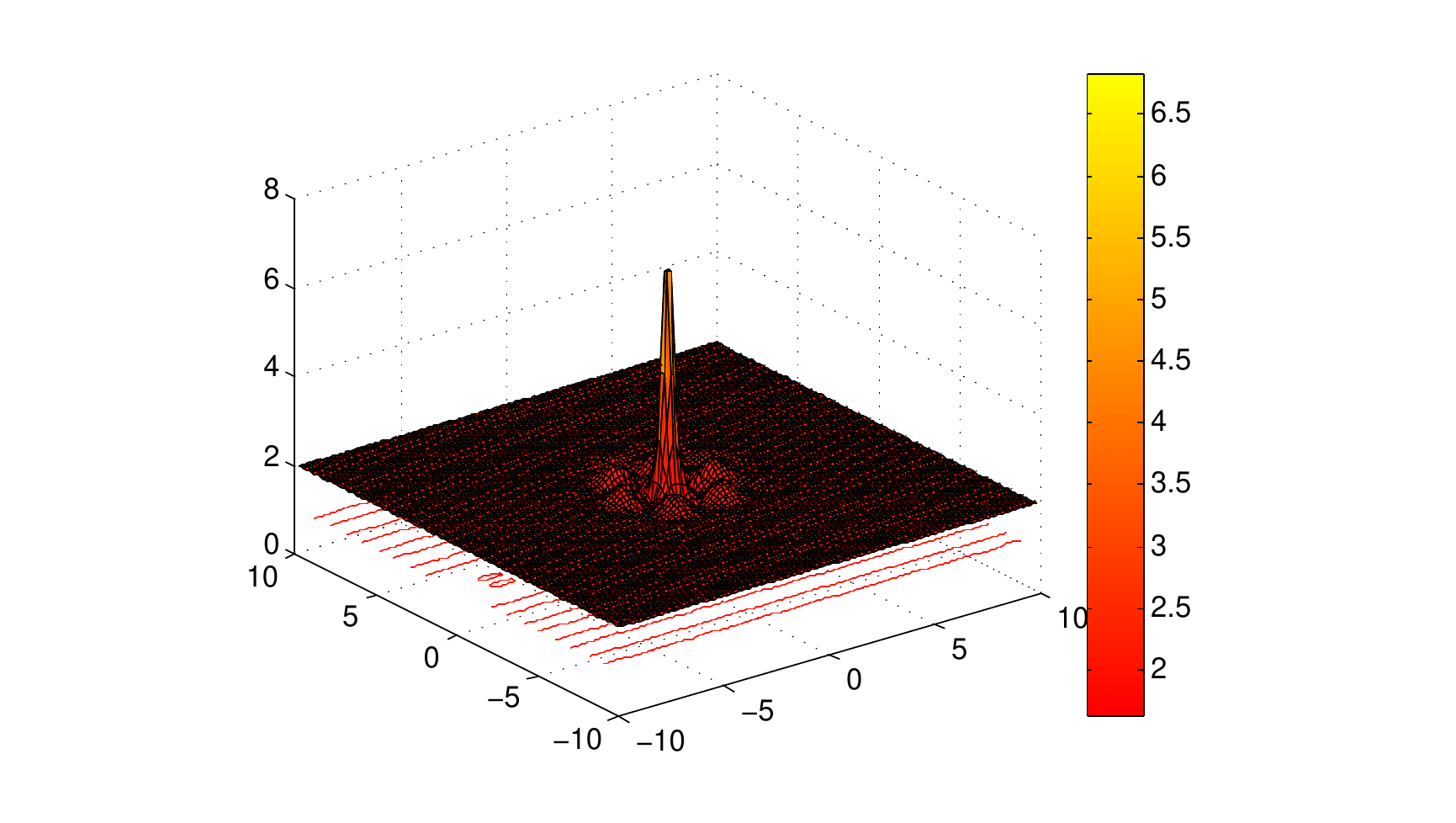} & 
\includegraphics[clip=true, trim= 2cm 9cm 2cm 9cm,  width=8cm]{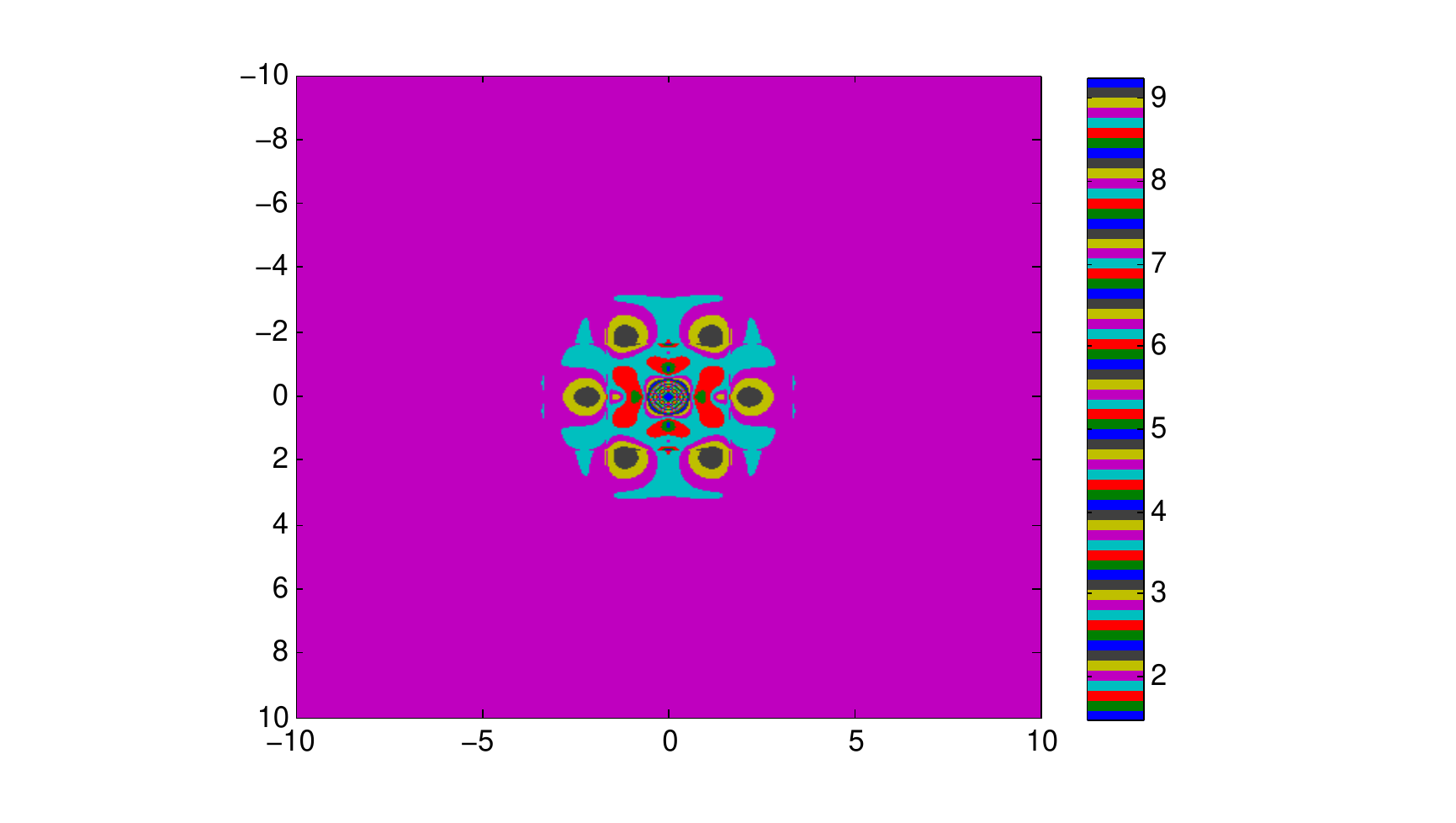} 
\end{array}
\end{displaymath}
\caption{Solution at $t=0.48$, with $v=10$.}
\end{figure}
\begin{figure}[p]
\begin{displaymath}
\begin{array}{cc}
\includegraphics[clip=true, trim= 2cm 9cm 2cm 9cm,  width=8cm]{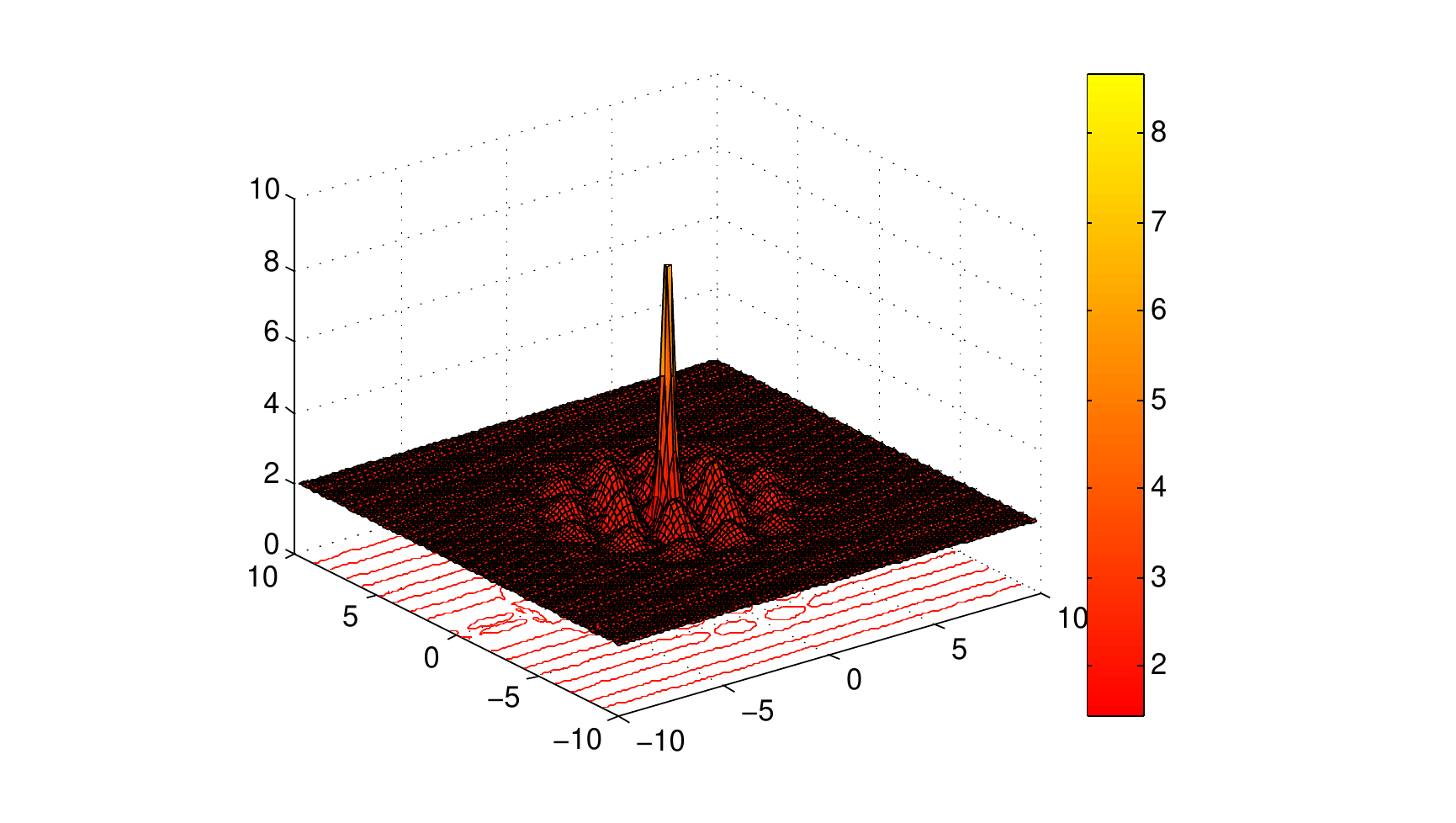} & 
\includegraphics[clip=true, trim= 2cm 9cm 2cm 9cm,  width=8cm]{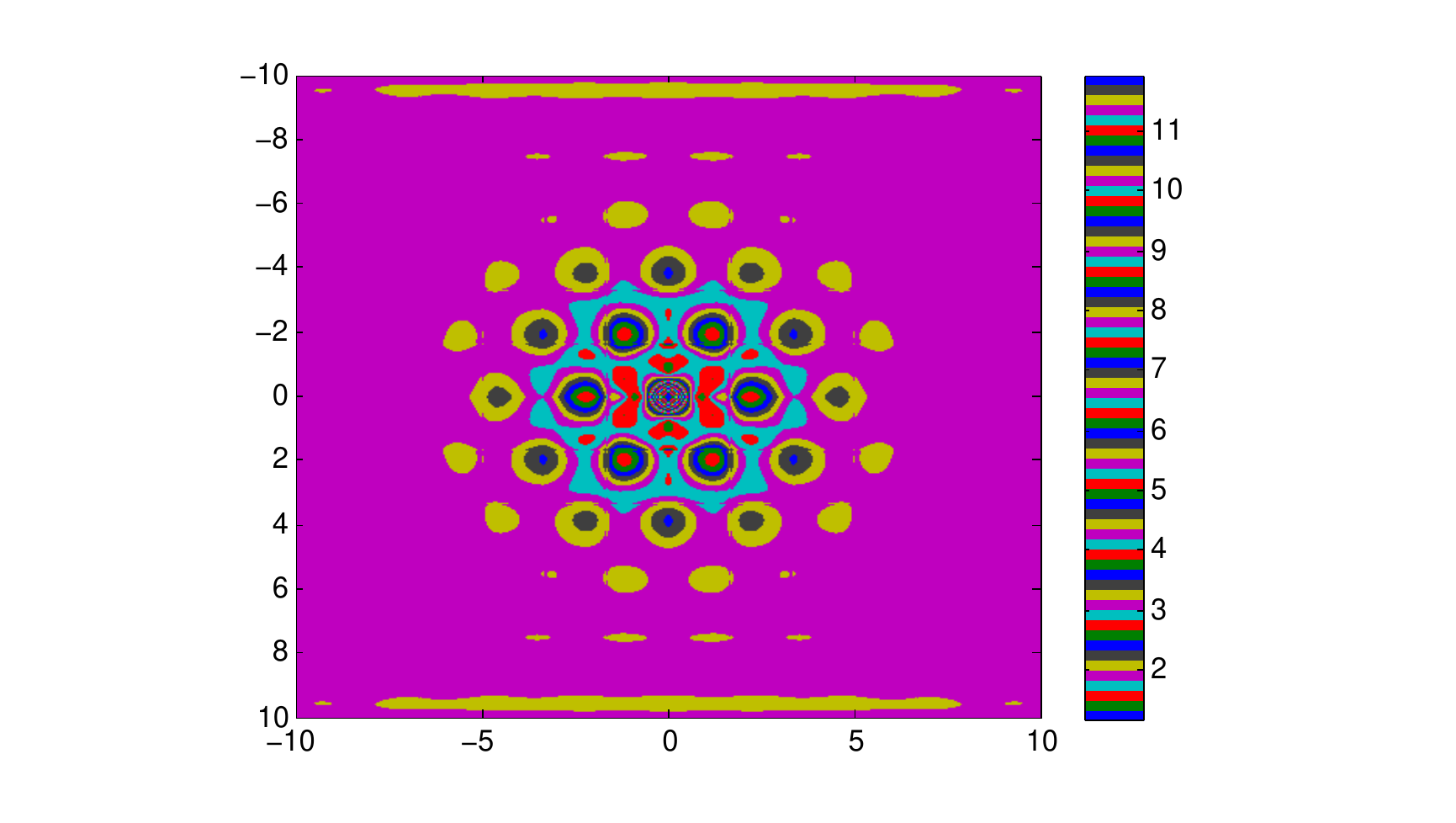} 
\end{array}
\end{displaymath}
\caption{Solution at $t=0.72$, with $v=10$.}
\end{figure}
\begin{figure}[p]
\begin{displaymath}
\begin{array}{cc}
\includegraphics[clip=true, trim= 2cm 9cm 2cm 9cm,  width=8cm]{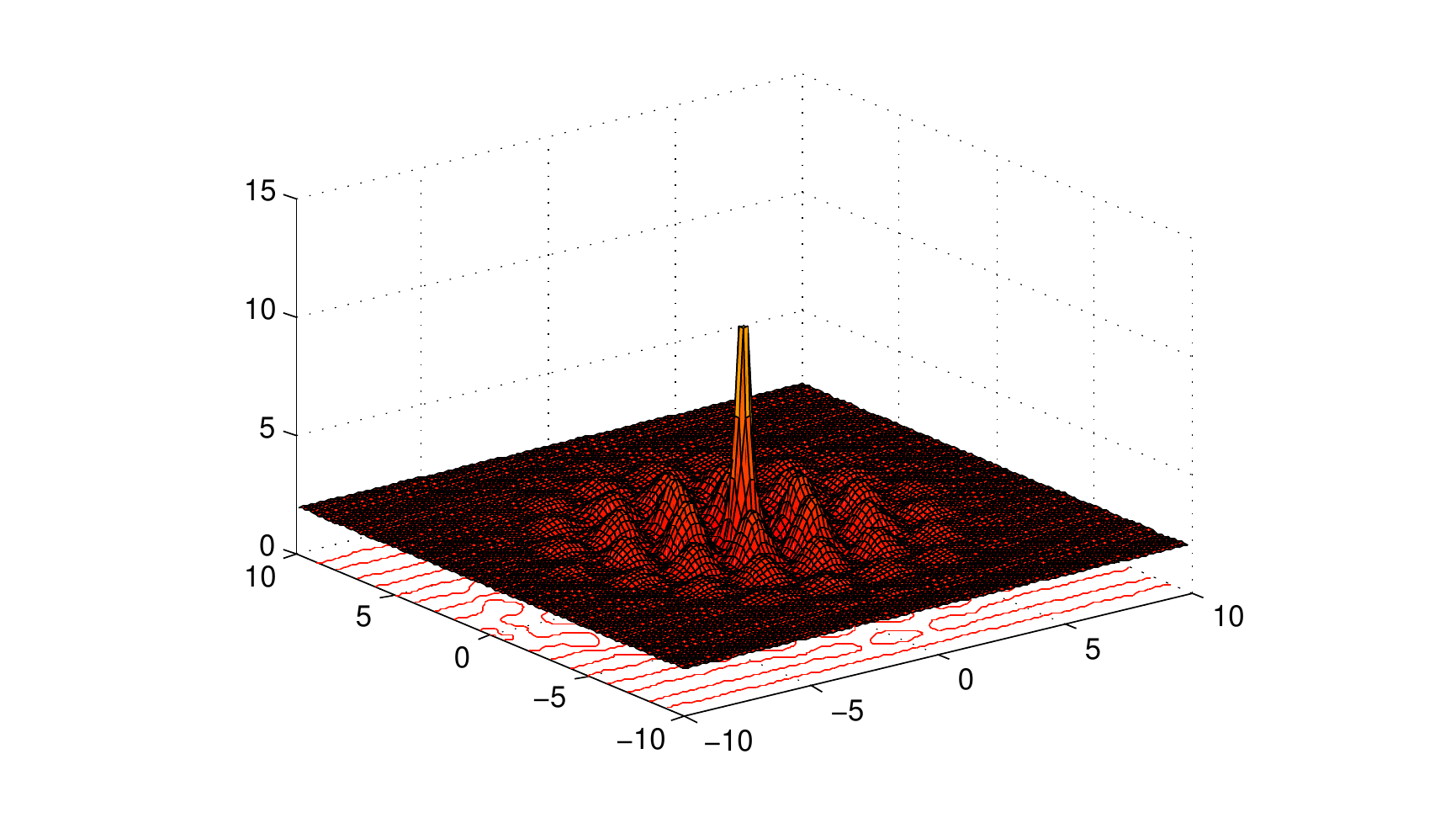} & 
\includegraphics[clip=true, trim= 2cm 9cm 2cm 9cm,  width=8cm]{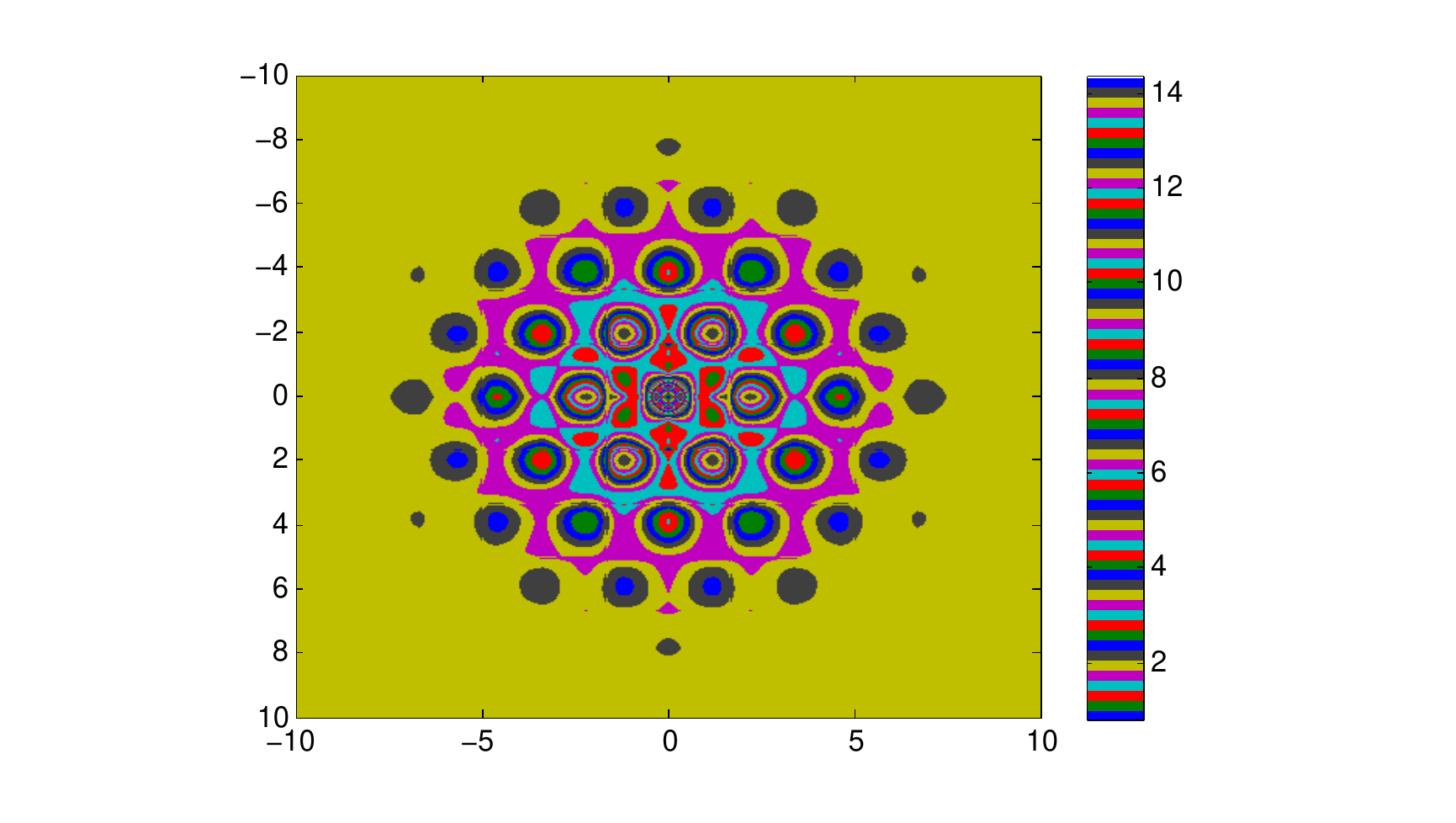} 
\end{array}
\end{displaymath}
\caption{Solution at $t=0.96$, with $v=10$.}
\end{figure}
\begin{figure}[p]
\begin{displaymath}
\begin{array}{cc}
\includegraphics[clip=true, trim= 2cm 9cm 2cm 9cm,  width=8cm]{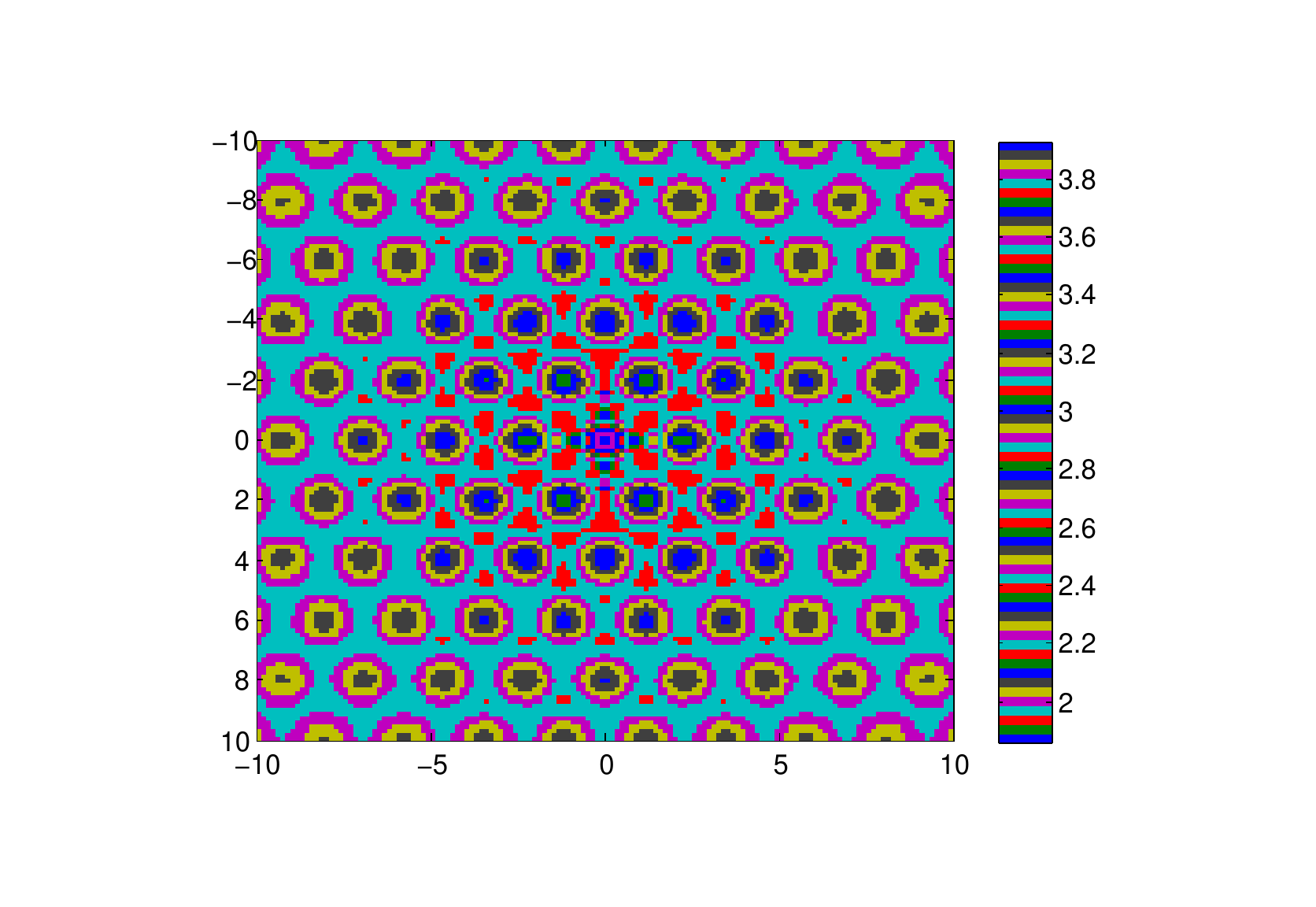} & 
\includegraphics[clip=true, trim= 2cm 9cm 2cm 9cm,  width=8cm]{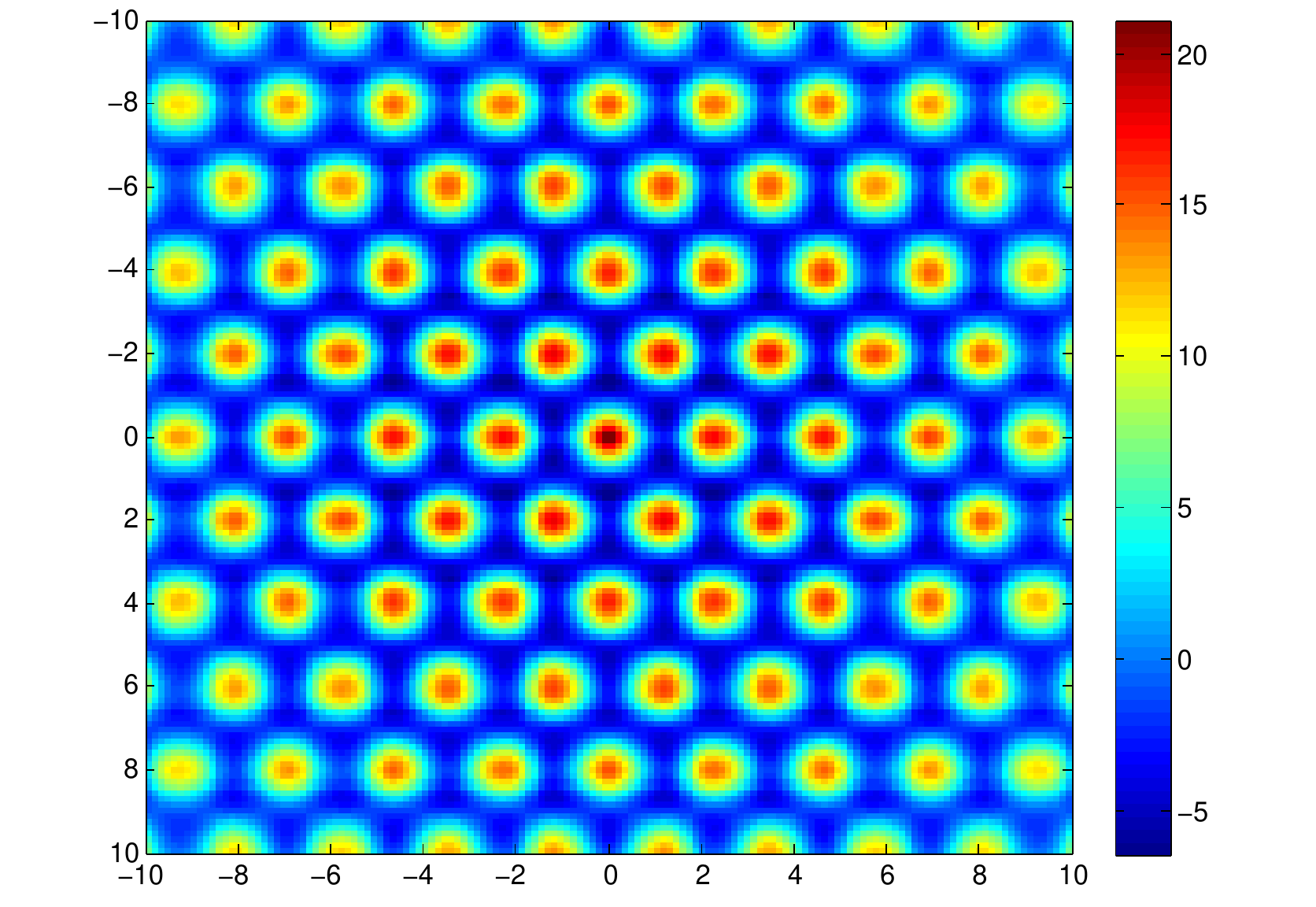} 
\end{array}
\end{displaymath}
\caption{Solution in the case of no delay, for $t=0.16$  (left-hand side) and $t=0.32$ (right-hand side).}
\end{figure}
\begin{table}
\begin{displaymath}
\begin{array} {|c|c|c|c|c|c|c|c|c|c|c|c|c|}
\hline
t & 0.08 & 0.16 & 0.24 & 0.32 & 0.4 & 0.48 & 0.56 & 0.64 & 0.72 & 0.80 & 0.88 & 0.96 \\ 
\hline
v=10 & 0.501 & 0.589 & 0.621 & 0.694 & 2.550  & 2.856 & 4.049 & 4.863 & 5.141  & 6.313 & 6.946 & 7.360 \\ 
v=20 & 0.501 & 0.589 & 0.621 & 2.812 & 4.724  & 6.207 & 7.224 & 8.625 & 10.455 & 11.577 & 12.93 & 13.11 \\
\makebox{no delay} & 0.501 & 6.207 & 14.14 & 14.14 & 14.14 & 14.14 & 14.14 & 14.14 &  14.14& 14.14 &
14.14 & 14.14 \\
\hline 
\end{array}
\end{displaymath}
\caption{ Dependence of the activity radius $r_a$ on $t$, in three different cases: $v=10$,$v=20$ and no delay.}
\end{table}
	
	{\bf Example 3.} In this example we consider a kind of neural field with spatially localized oscillations which occur in excitable neural media with short-range 
		excitation and long-range inhibition (mexican hat connectivity), in the case of a spatially localized input.
		Such neural fields are known as {\em breathers} and their dynamics are analysed in \cite{Bressloff}. A similar example was described and computed
		in\cite{Hutt2010}. In this case the
	 connectivity kernel has the form: 
	\begin{equation} \label{ker3}
  K(r)= 20	\frac{\exp (-r)}{18 \pi}-14 \frac{\exp \left(- \frac{-r}{3} \right)}{18 \pi}.
	\end{equation}
The firing rate function has the  the sigmoidal form
\[
S(x)= \frac{2}{1+\exp(-10000(x-0.005))}
\]
and the external input is given by
\[
I(\bar{x},t)= 5 \frac{\exp (-x^2/32 -y^2/32)}{32 \pi}.
\]
Our aim is to investigate how the dynamic behaviour of the solution is influenced by the propagation speed of the impulses, so
we consider different values of $v$. As initial condition we  consider $V_0 \equiv 0$. Note that when there is no delay (infinite propagation speed) the solution tends to a stable stationary state.
In the case of finite propagation speed, oscillatory solutions occur. The frequency of the oscillations
decreases as the propagation speed becomes lower. As an example, in Fig.22 the graphs of the 
solution maximum and minimum are depicted, as a function  of time, in the case $v=50$.
We have carried out  simulations over the time interval $[0,1.6]$ with stepsize in time $h_t=0.02$.
The parameters of the space discretization were $m=24$, $N=48$.

\begin{figure}[p]
\includegraphics[clip=true, trim= 2cm 9cm 2cm 9cm,  width=14cm]{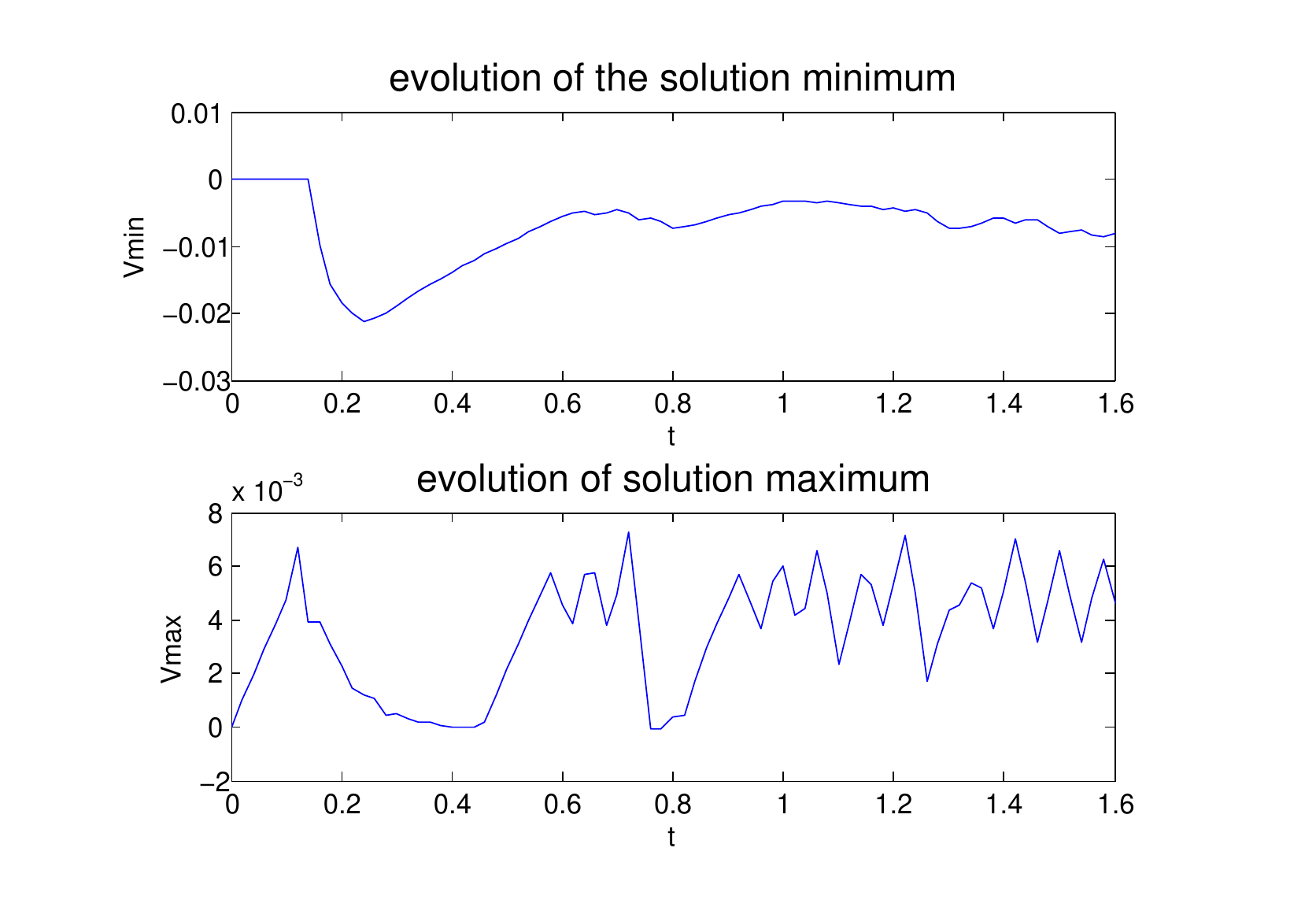} 
\caption{Evolution of the solution maximum and minimum for the breather with $v=50$.}
\end{figure}
In Fig. 23 and 24, the 3d-graphs of the solution, in the case $v=50$, are displayed for different
moments of time. From the pictures it is clear that the solution has oscillatory behaviour and the amplitude of oscillations 
decreases with time. As referred in \cite{Bressloff}, we also can see a symmetry breaking dynamical instability, when a two-dimensional
radially symmetric stationary pulse bifurcates to a nonradially symmetric breather 
with an integer number of lobes (in this case 4). 
%
\begin{figure}[p]
\includegraphics[clip=true, trim= 2cm 2cm 2cm 2cm,  width=12cm]{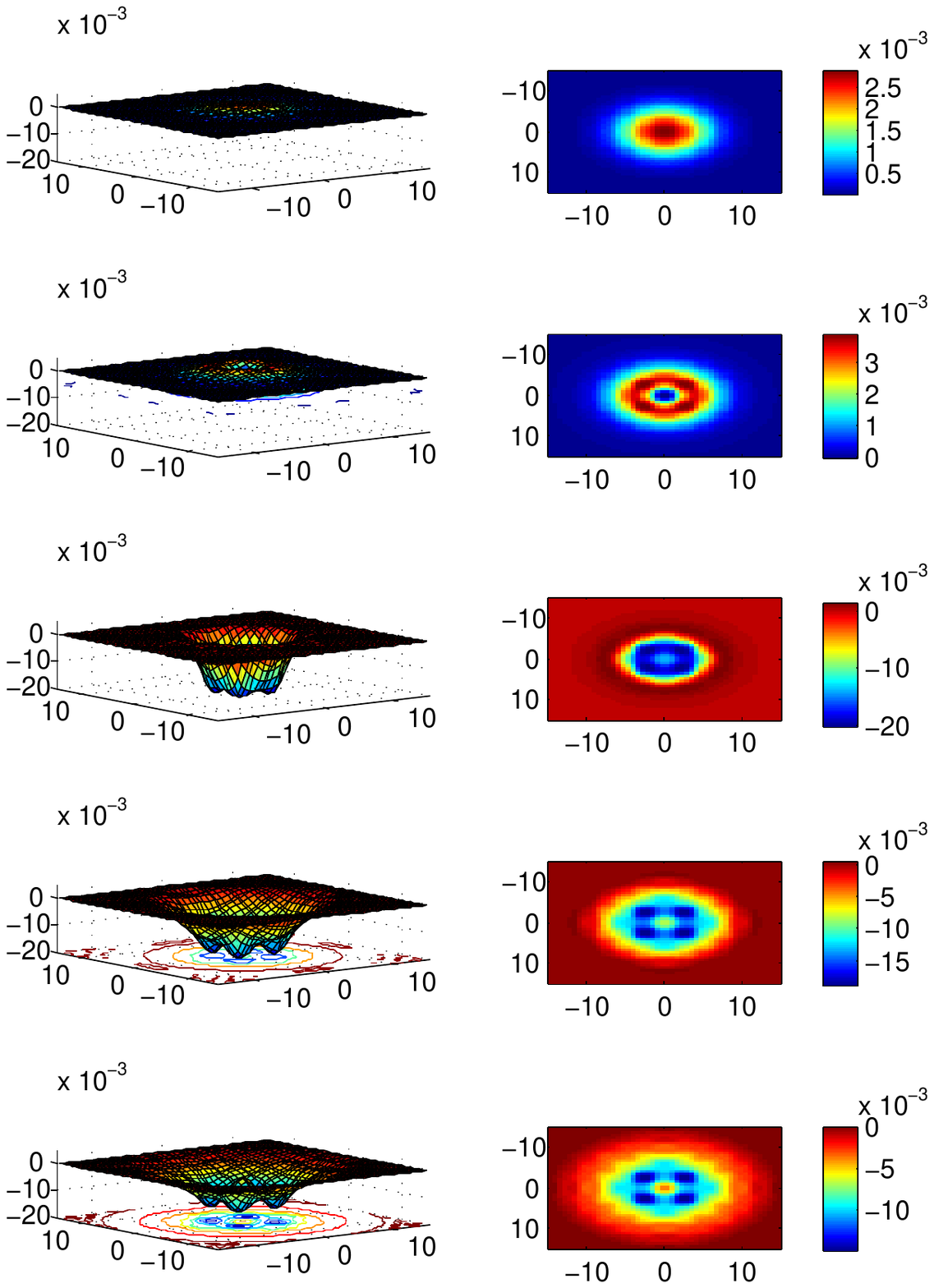} 
\caption{Plots of the solution in the case $v=50$ at different moments of time: $t=0.08,0.16,0.24, 0.32, 0.40$ .}
\end{figure}
\begin{figure}[p]
\includegraphics[clip=true, trim= 2cm 2cm 2cm 2cm,  width=12cm]{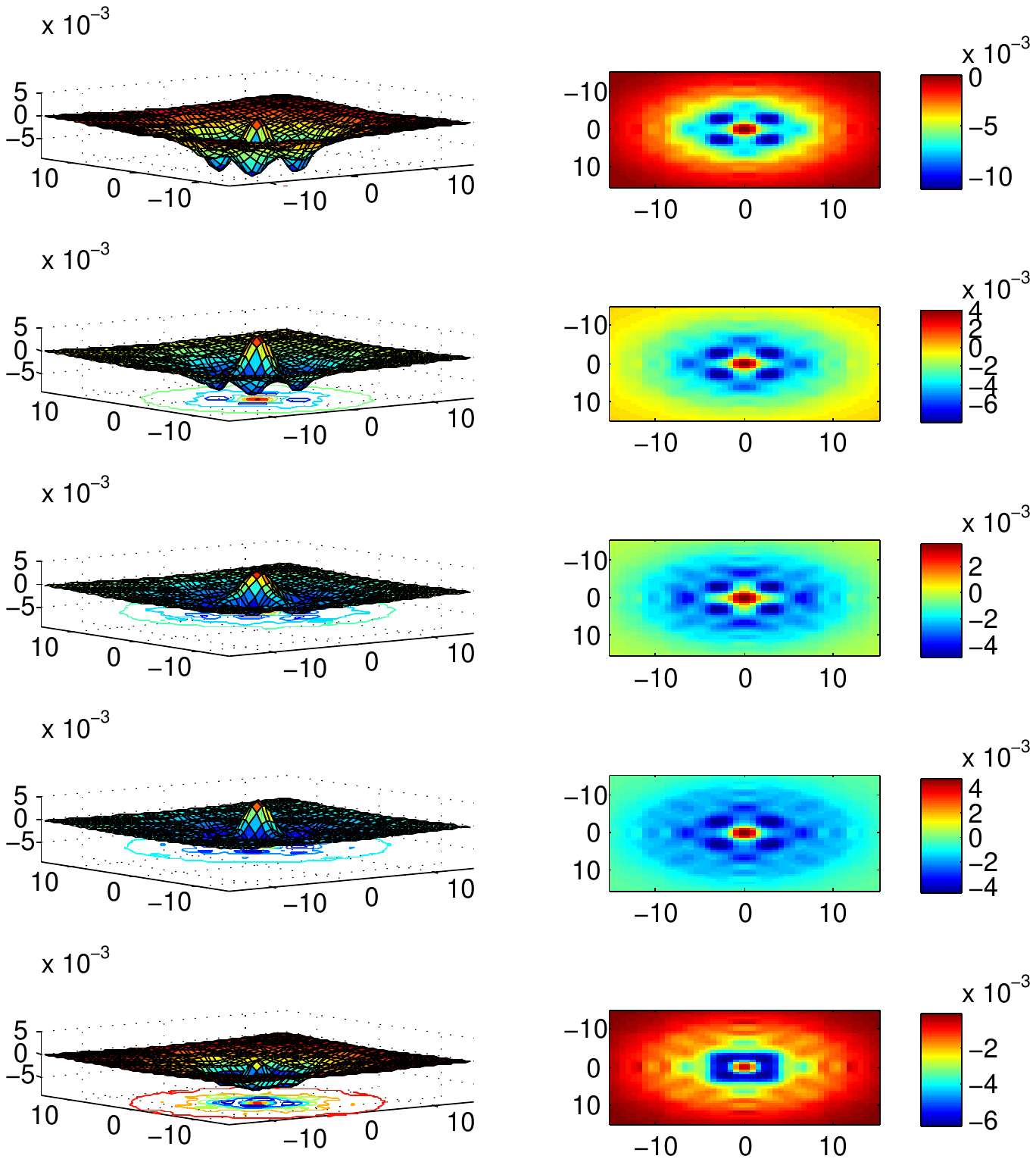} 
\caption{Plots of the solution in the case $v=50$ at different moments of time: $t=0.48,0.56,0.64, 0.72, 
0.80$ .}
\end{figure}
\section{Conclusions}
In this paper we present some numerical simulations obtained by an algorithm which was previously described in \cite{LP1} and \cite{LP2}.
The three numerical examples considered here have been analysed before by other authors 
and are directly connected with real world applications.
The main purposes are 1) to test the performance of the mentioned algorithm, by comparing the numerical results with those obtained by other authors;  2) to analyse  with more detail the properties of the solutions and take conclusions about their physical meaning.
The numerical experiments have shown that the computational method applied in these simulations is an efficient tool and provides sufficiently accurate results to analyse the physical phenomena under consideration.
All the numerical simulations were carried out in a personal computer with a Intel Processor of 2.4 Ghz.
The computing time for each time step is about 10-20s, depending on the example considered.

The interpolation method used to reduce the matrix size has proved to be very efficient in the cases where
the input data are sufficiently smooth. However, in the case of an external input which changes very fast in the neighborhood of a certain point (like in example~2 ) this method may not be applicable.
In such cases the algorithm can still be used but the computing time is much higher.
This suggests  that it would be useful in such cases to consider the use of nonuniform meshes, which take into account the specific behaviour of the solution. Such improvement of the algorithm will be the subject of future research.

\section{Acknowledgements}
This research was supported by a Marie Curie Intra European Fellowship within the 7th European Community Framework Programme (PIEF-GA-2013-629496).


\begin{thebibliography}{99}
\bibitem{Amari} S.L. Amari, Dynamics of pattern formation in lateral-inhibition type
neural fields, Biol. Cybernet. 27 (2) (1977) 77--87.

\bibitem{Bressloff} S.E. Folias and P.C.Bressloff, Breathers in Two-dimensional neural media, Phys. Rev. Let. 95 (2005) 208107. 

\bibitem{Brunner} H. Brunner, On the numerical solution of nonlinear Volterra-Fredholm
integral equations by collocation methods, SIAM J. Numer. Anal. 27 (1990) 987-1000.
\bibitem{Cardone} A. Cardone, E. Messina, and E.Russo, A fast iterative method for discretized Volterra-Fredholm integral equations, J. Comput. Applied Math. 189 (2006)
568-579.
\bibitem{Faye} G. Faye and O. Faugeras, Some theoretical and numerical results for 
delayed neural field equations, Physica D 239 (2010) 561--578.
\bibitem{Han} Han Guoqiang, Asymtotic error expansion for a nonlinear
Volterra-Fredholm integral equation, J. Comput. Applied Math. 59 (1995) 49-59.

\bibitem{Hutt2010} A. Hutt and N. Rougier, Activity spread and breathers
induced by finite transmission speeds in two-dimensional neuronal fields,
Physical Review,E 82 (2010) 055701.

\bibitem{Hutt2014} A. Hutt and N. Rougier, Numerical Simulations of One-
and Two-dimensional Neural Fields Involving Space-Dependent Delays,
in S. Coombes et al., Eds., Neural Fields Theory and Applications,
Springer, 2014.

\bibitem{Jacki1} Z. Jackiewicz, M. Rahman, B.D. Welfert, Numerical Solution of
a Fredholm integro-differential equation modelling neural networks, 
Applied Numerical Mathematics, 56 (2006) 423-432.

\bibitem{Jacki2} Z. Jackiewicz, M. Rahman, B.D. Welfert, Numerical Solution of
a Fredholm integro-differential equation modelling $\theta $-neural networks, 
Appl. Math. Comput., 195 (2008) 523-536.

\bibitem{Kauthen} J.-P. Kauthen, Continuous time collocation methods for Volterra-Fredholm integral equations, Num. Math. {\bf 56} (1989) 409--424.

\bibitem{LP1} P. M. Lima and E. Buckwar, Numerical solution of the neural field equation in the two-dimensional case, Arxiv ref.  http://arxiv.org/abs/1508.07484, 2015.
\bibitem{LP2} P.M. Lima and E. Buckwar, Numerical investigation of the two-dimensional  neural field equation with delay,  to appear in Proceedings of IEEE.

\bibitem{Lund} J. Lund, A. Angelucci, P. Bressloff, Anatomical substracts for functional columns in macaque 
monkey primary visual cortex,  Cerebral Cortex, {\bf 13} (2003) 15-24.
\bibitem {Pothast} R. Pothast and P. beim Graben, Existence and properties of solutions
for neural field equations, Math. Meth. Appl. Sci., 33 (2010) 935-949.


\bibitem{WC} H.R. Wilson and J.D. Cowan, Excitatory and inhibitory interactions
in localized populations of model neurons, Bipophys. J., 12 (1972) 1-24.


\bibitem{Xie} Weng-Jing Xie, Fu-Rong Lin, A fast numerical solution method
for two dimensional Fredholm integral equations of the second kind,
Applied Numerical Mathematics, 59 (2009) 1709-1719.

\end{thebibliography}
\end{document}